\RequirePackage{fix-cm}
\documentclass[smallcondensed,final]{svjour3}     

\usepackage{moreverb}

\usepackage[colorlinks,bookmarksopen,bookmarksnumbered,citecolor=red,urlcolor=red]{hyperref}

\usepackage{algpseudocode} 
\usepackage[section]{algorithm}

\usepackage{tikz}
\usetikzlibrary{shapes,arrows,trees,snakes} 
\usepackage{pgfplots}
\pgfplotsset{compat=1.3}
\usepgfplotslibrary{external} 
\tikzset{external/system call={lualatex \tikzexternalcheckshellescape -halt-on-error -interaction=batchmode -jobname "\image" "\texsource"}}
\usepackage{amsmath,amsfonts,amssymb}
\usepackage{cite}
\usepackage{mathtools}

\definecolor{utorange}{RGB}{203,96,21}
\definecolor{utblack}{RGB}{99,102,106}
\definecolor{utbrown}{RGB}{110,98,89}
\definecolor{utsecbrown}{RGB}{217,200,158}
\definecolor{utsecgreen}{RGB}{208,222,187}
\definecolor{utsecblue}{RGB}{127,169,174}

\newcommand{\bs}[1]{\ensuremath{\boldsymbol #1}}

\newcommand\BibTeX{{\rmfamily B\kern-.05em \textsc{i\kern-.025em b}\kern-.08em
T\kern-.1667em\lower.7ex\hbox{E}\kern-.125emX}}

\begin{document}

\title{Comparison of Multigrid Algorithms for High-order
  Continuous Finite Element Discretizations}
\titlerunning{Comparison of multigrid algorithms for high-order discretizations}
\author{Hari Sundar \and Georg Stadler \and George Biros}
\authorrunning{Hari Sundar, Georg Stadler and George Biros}

\institute{Hari Sundar \at
 School of Computing, University of Utah, Salt Lake City, UT\\
\email{hari@cs.utah.edu}
\and Georg Stadler \at
Courant Institute of Mathematical Sciences, New York University, New York, NY
\and George Biros\at
Institute for Computational Engineering and Sciences, The University of Texas at Austin, Austin, TX
}
\date{}

\maketitle

\vspace{-0.5in}
\begin{abstract}
We present a comparison of different multigrid approaches for the
solution of systems arising from high-order continuous finite element
discretizations of elliptic partial differential equations on complex
geometries.  We consider the pointwise Jacobi, the
Chebyshev-accelerated Jacobi and the symmetric successive
over-relaxation (SSOR) smoothers, as well as elementwise block Jacobi
smoothing. Three approaches for the multigrid hierarchy are compared:
(1) high-order $h$-multigrid, which uses high-order interpolation and
restriction between geometrically coarsened meshes; (2) $p$-multigrid,
in which the polynomial order is reduced while the mesh remains
unchanged, and the interpolation and restriction incorporate the
different-order basis functions; and (3), a first-order approximation
multigrid preconditioner constructed using the nodes of the high-order
discretization.  This latter approach is often combined with algebraic
multigrid for the low-order operator and is attractive for high-order
discretizations on unstructured meshes, where geometric coarsening is
difficult.  Based on a simple performance model, we compare the
computational cost of the different approaches.  Using scalar test
problems in two and three dimensions with constant and varying
coefficients, we compare the performance of the different multigrid
approaches for polynomial orders up to 16. Overall, both $h$- and
$p$-multigrid work well; the first-order approximation is less
efficient. For constant coefficients, all smoothers work well. For
variable coefficients, Chebyshev and SSOR smoothing outperform Jacobi
smoothing.  While all of the tested methods converge in a
mesh-independent number of iterations, none of them behaves completely
independent of the polynomial order.  When multigrid is used as a
preconditioner in a Krylov method, the iteration number decreases
significantly
compared to using multigrid as a solver.
\end{abstract}

\section{Introduction}

This paper presents a comparison of geometric multigrid
methods for the solution of systems arising from high-order (we target
polynomial orders up to 16) continuous finite element discretizations
of elliptic partial differential equations. Our particular interest is
to compare the efficiency of different multigrid methods for elliptic
problems with varying coefficients on complex geometries.
High-order spatial discretizations for these problems can have significant advantages
over low-order methods since they reduce the problem size for given
accuracy, and allow for better performance on modern hardware.
The main challenges in high-order discretizations are that matrices
are denser compared to low-order methods, and that they lose structural
properties such as the M-matrix
property, which often allows to prove convergence of iterative
solvers.

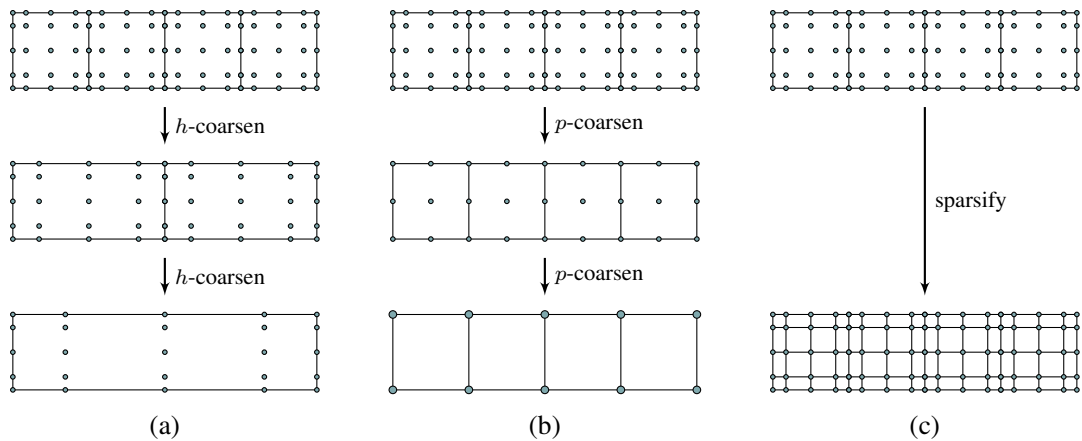
\begin{figure}[t]
		\begin{tikzpicture}[scale=1.0]
		\draw (-5,4) grid +(4,1);
		\foreach \e in {-5,...,-2}
		\foreach \x in {0,0.1727,0.5,0.8273, 1.0} {
			\draw[fill=utsecblue] (\e+\x, 4) circle (0.03);
			\draw[fill=utsecblue] (\e+\x, 4.1727) circle (0.03);
			\draw[fill=utsecblue] (\e+\x, 4.5) circle (0.03);
			\draw[fill=utsecblue] (\e+\x, 4.8273) circle (0.03);
			\draw[fill=utsecblue] (\e+\x, 5) circle (0.03);
		}
		\draw[-latex',thick] (-3, 3.75) -- node[right] {{\scriptsize $h$-coarsen}} (-3, 3.25);
		\draw (-5,2) rectangle +(4,1);
		\draw (-3,2) -- (-3,3);
		\foreach \e in {-5,-3}
		\foreach \x in {0,0.1727,0.5,0.8273, 1.0} {
			\draw[fill=utsecblue] (\e+2*\x, 2) circle (0.03);
			\draw[fill=utsecblue] (\e+2*\x, 2.1727) circle (0.03);
			\draw[fill=utsecblue] (\e+2*\x, 2.5) circle (0.03);
			\draw[fill=utsecblue] (\e+2*\x, 2.8273) circle (0.03);
			\draw[fill=utsecblue] (\e+2*\x, 3) circle (0.03);
		}
	
		\draw[-latex',thick] (-3, 1.75) -- node[right] {{\scriptsize $h$-coarsen}} (-3, 1.25);
	
		\draw (-5,0) rectangle +(4,1);
		\foreach \x in {0,0.1727,0.5,0.8273, 1.0} {
			\draw[fill=utsecblue] (-5+4*\x, 0) circle (0.03);
			\draw[fill=utsecblue] (-5+4*\x, 0.1727) circle (0.03);
			\draw[fill=utsecblue] (-5+4*\x, 0.5) circle (0.03);
			\draw[fill=utsecblue] (-5+4*\x, 0.8273) circle (0.03);
			\draw[fill=utsecblue] (-5+4*\x, 1) circle (0.03);
		}
		
		\draw (0,4) grid +(4,1);
		\foreach \e in {0,...,3}
		\foreach \x in {0,0.1727,0.5,0.8273, 1.0} {
			\draw[fill=utsecblue] (\e+\x, 4) circle (0.03);
			\draw[fill=utsecblue] (\e+\x, 4.1727) circle (0.03);
			\draw[fill=utsecblue] (\e+\x, 4.5) circle (0.03);
			\draw[fill=utsecblue] (\e+\x, 4.8273) circle (0.03);
			\draw[fill=utsecblue] (\e+\x, 5) circle (0.03);
		}
	
		\draw[-latex',thick] (2, 3.75) -- node[right] {{\scriptsize $p$-coarsen}} (2, 3.25);
	
		\draw (0,2) grid +(4,1);
		\foreach \x in {0,0.5,...,4} {
			\draw[fill=utsecblue] (\x, 2) circle (0.03);
			\draw[fill=utsecblue] (\x, 2.5) circle (0.03);
			\draw[fill=utsecblue] (\x, 3) circle (0.03);
		}
	
		\draw[-latex',thick] (2, 1.75) -- node[right] {{\scriptsize $p$-coarsen}} (2, 1.25);
	
		\draw (0,0) grid +(4,1);
		\foreach \x in {0,1,2,3,4} {
			\draw[fill=utsecblue] (\x, 0) circle (0.05);
			\draw[fill=utsecblue] (\x, 1) circle (0.05);
		}
		
			\draw (5,4) grid +(4,1);
			\foreach \e in {5,...,8}
			\foreach \x in {0,0.1727,0.5,0.8273, 1.0} {
				\draw[fill=utsecblue] (\e+\x, 4) circle (0.03);
				\draw[fill=utsecblue] (\e+\x, 4.1727) circle (0.03);
				\draw[fill=utsecblue] (\e+\x, 4.5) circle (0.03);
				\draw[fill=utsecblue] (\e+\x, 4.8273) circle (0.03);
				\draw[fill=utsecblue] (\e+\x, 5) circle (0.03);
			}
	
			\draw[-latex',thick] (7, 3.75) -- node[right] {{\scriptsize sparsify}} (7, 1.25);
	
			\draw[step=0.5] (4.99,0) grid +(4.01,1);
			\draw (5,0.1727) -- (9,0.1727);
			\draw (5,0.8273) -- (9,0.8273);
			\foreach \e in {5,...,8} {
				\draw (\e+0.1727,0) -- (\e+0.1727,1);
				\draw (\e+0.8273,0) -- (\e+0.8273,1);
				\foreach \x in {0,0.1727,0.5,0.8273, 1.0} {
					\draw[fill=utsecblue] (\e+\x, 0) circle (0.03);
					\draw[fill=utsecblue] (\e+\x, 0.1727) circle (0.03);
					\draw[fill=utsecblue] (\e+\x, 0.5) circle (0.03);
					\draw[fill=utsecblue] (\e+\x, 0.8273) circle (0.03);
					\draw[fill=utsecblue] (\e+\x, 1) circle (0.03);
				}
			}
      \node at (-3, -0.5) { (a) };
      \node at (2, -0.5) { (b) };
      \node at (7, -0.5) { (c) };
		
    \end{tikzpicture}
		\caption{\label{fig:approaches} Illustration of
                  different multigrid hierarchies for high-order
                  finite element discretizations: (a) high-order
                  $h$-multigrid, (b) $p$-multigrid and (c)
                  low-order approximation preconditioner based
                  on the nodes of the high-order discretization
                  .}
\end{figure}

As illustrated in Figure~\ref{fig:approaches}, there are several
possibilities for constructing a multigrid hierarchy for high-order
discretizations: (1) high-order geometric $h$-multigrid, where the
mesh is coarsened geometrically and high-order interpolation and
prolongation operators are used; (2) $p$-multigrid, in which the
problem is coarsened by reducing the polynomial order, and the
interpolation and prolongation take into account the different order
basis functions; and (3) a first-order approximation as
preconditioner, constructed from the nodes of the high-order
discretization.
For the polynomial orders $1\le p\le 16$, we
compare these multigrid approaches, combined with different
smoothers. We also compare the use of multigrid as a solver as well as a preconditioner
in a Krylov subspace method.  While we use moderate size model
problems (up to about $2$~million unknowns in 3D), we also discuss our findings with regard to parallel
implementations on high performance computing platforms. 
We also discuss parallelization aspects relevant for implementations on
shared or distributed memory architectures. For instance, the
implementation of Gauss-Seidel smoothers can be challenging in
parallel~\cite{AdamsBrezinaHuEtAl03, BakerFalgoutKolevEtAl11}; for
this reason, we include a Chebyshev-accelerated Jacobi smoother in our
comparisons. This Chebyshev smoother is easy to implement in parallel,
and often is as effective a smoother as Gauss-Seidel.


We use high-order discretizations based on Legende-Gauss-Lobotto
(LGL) nodal basis functions on quadrilateral or hexahedral
meshes. Tensorized basis functions allow for a fast, matrix-free
application of element matrices. This is particularly important for
high polynomial degrees in three dimensions, as element matrices
can become large. For instance, for a three-dimensional hexahedral mesh
and finite element discretizations with polynomial degree $p$, the
dense element
matrices are of size $(p+1)^3\times (p+1)^3$. Thus, for $p=8$, this
amounts to more than half a million entries per element.  For
tensorized nodal basis functions on hexahedral meshes, the application
of elemental matrices to vectors can be implemented efficiently by
exploiting the tensor structure of the basis functions, as is common
for spectral elements ~\cite{DevilleFischerMund02}.

{\em Related work:} Multigrid for high-order/spectral finite elements
has been studied as early as in the 1980s. In~\cite{RonquistPatera87},
the authors observe that point smoothers such as the simple Jacobi
method result in resolution-independent convergence rates for
high-order elements on simple one and two-dimensional
geometries. Initial theoretical evidence for this behavior is given
in~\cite{MadayMunoz88}, where multigrid convergence is studied for
one-dimensional spectral methods and spectral element problems.
The use of $p$-multigrid is rather common in the context of
high-order discontinuous Galerkin discretizations
\cite{FidkowskiOliverLuEtAl05, HelenbrookAtkins06}, but $p$-multigrid
has also been used for continuous finite element discretizations ~\cite{HelenbrookMavriplisAtkins03}.
A popular strategy for high-order discretizations on unstructured
meshes, for which geometric mesh coarsening is
challenging, is to assemble a low-order approximation of the
high-order system and use an algebraic multigrid method to invert the
low-order (and thus much sparser) operator~\cite{Brown10, Kim07,
  DevilleMund90, Olson07, CanutoGervasioQuarteroni10}.
In~\cite{HeysManteuffelMcCormickEtAl05}, this approach is compared
with the direct application of algebraic multigrid to the high-order
operator and the authors find that one of the main difficulties is the
assembly of the high-order matrices required by algebraic multigrid
methods.


{\em Contributions:} 
There has been a lot of work on high-order discretization methods and on the efficient application of the resulting operators. However, efficient solvers for such
discretization schemes have received much less attention. In particular, theoretical and experimental studies are scattered regarding the actual performance (say the number of v-cycles or matrix-vector products to solve a system) of
the different schemes under different scenarios. A systematic analysis of such performance is not available.
In this paper, we address this gap in the existing literature. In particular we (1) consider high-order continuous Galerkin discretizations up to 16th order, (2) examine three different multigrid hierarchies ($h$, $p$, and first-order), (3) examine several different smoothers: Jacobi, polynomial, SSOR, and block Jacobi, (4) consider different settings (constant, mildly variable, and highly variable) of coefficients and (5) consider problems in 2D and 3D.
To our knowledge, this is the first study of this kind. Our results demonstrate significant variability in the performance of the different schemes for higher-order
elements, highlighting the need for further research on the smoothers. 
Although the overall runtime will depend on several factors---including the implementation and the target architecture---in this work we limit ourselves to characterizing performance as the number of fine-grid matrix-vector products needed for convergence. This is the most dominant cost and is also independent of the implementaion and architecture, allowing for easier interpretation and systematic comparison with other approaches. 
Finally, we provide an easily extendable Matlab
implementation,\footnote{\url{http://hsundar.github.io/homg/}} which allows a systematic comparison of the different methods in the same framework.


{\em Limitations:} While this work is partly driven by our interest in
scalable parallel simulations on nonconforming meshes derived from
adaptive octrees (e.g.,\cite{SundarBirosBursteddeEtAl12}), 
 for the comparisons
presented in this paper we restrict ourselves to moderate size
problems on conforming meshes. We do not fully address time-to-solution, as
we do not use a high-performance implementation. However, recent
results using a scalable parallel implementation indicate that many of
our observations generalize to non-conforming meshes and that the
methods are scalable to large parallel computers
\cite{GholaminejadMalhotraSundarEtAl14}. While we
experiment with problems with strongly varying coefficients, we do not
study problems with discontinuous or anisotropic coefficients, nor consider
ill-shaped elements.

{\em Organization of this paper:} In \S\ref{sec:problem} we describe
the test problem, as well as discretization approach for the different
multigrid schemes. In \S\ref{sec:approaches}, we describe in detail
the different multilevel approaches for solving the resulting
high-order systems. In \S\ref{sec:numerics}, we present a
comprehensive comparison of different approaches using test problems
in 2D and 3D. Finally, in \S\ref{sec:discuss} we draw conclusions and
discuss our findings.

\section{Problem statement and preliminaries}
\label{sec:problem}

We wish to solve the Poisson problem with
homogeneous Dirichlet boundary conditions on an open bounded domain
$\Omega\subset\mathbb R^d$ ($d=2$ or $d=3$) with boundary $\partial
\Omega$, i.e., we search the solution $u(\bs x)$ of:
\begin{equation}\label{eq:Poisson}
  \begin{aligned}
    -\nabla\cdot\left(\mu(\bs x)\nabla u(\bs x)\right) &= f(\bs x) \quad &&\text{ for } \bs x\in \Omega,\\
    \quad u(\bs x)& = 0  \quad &&\text{ for } \bs x\in \partial\Omega.
  \end{aligned}
\end{equation}
Here, $\mu(\bs x)\ge \mu_0>0$ is a spatially varying coefficient that
is bounded away from zero, and $f(x)$ is a given right hand side. We
discretize~\eqref{eq:Poisson} using finite elements with basis
functions of polynomial order $p$ and solve the resulting discrete
system using different multigrid variants. Next, in
\S\ref{subsec:galerkin} and \S\ref{sub:restriction_&_prolongation}, we
discuss the Galerkin approximation to~\eqref{eq:Poisson} and the setup
of the inter-grid transfer operators to establish a multilevel
hierarchy. In \S\ref{sub:meshing}, we discuss details of the meshes
and implementation used for our comparisons.

\subsection{Galerkin approximation} \label{subsec:galerkin}

Given a bounded, symmetric bilinear form\footnote{In our case,
$a(u,v)=\int_\Omega \mu\nabla u \cdot \nabla v$.} $a(u,v)$ that is
coercive
on $H_0^{1}(\Omega)$, and $f \in L^{2}(\Omega)$, we want to find $u
\in H_0^{1}(\Omega)$ such that $u$ satisfies
\begin{equation}
\label{eqn:weakForm}
a(u,v) =  (f,v)_{L^2(\Omega)}, \ \ \ \forall v \in H_0^{1}(\Omega),
\end{equation}
where $(f,v)_{L^2(\Omega)} = \int_\Omega fv\,dx$ and
$H_0^1(\Omega)\subset L^2(\Omega)$ denotes the subspace of
functions with square integrable derivatives that vanish on the
boundary.
This problem is known to have a unique solution $u^*$ \cite{BrennerScott94}. 
We now derive discrete equations whose solutions
approximate the solution of
(\ref{eqn:weakForm}). First, we define a sequence of $m$ nested conforming
{\em finite}-dimensional spaces, $V_1 \subset V_2 \subset \cdots \subset V_m \subset
H_0^{1}(\Omega)$.
Here, $V_k$ is the finite element space that corresponds to a finite element mesh
at a specified
polynomial order, and $V_{k-1}$ corresponds to the next coarser
problem,
as illustrated  in Figure~\ref{fig:approaches}-(a,b) for different coarsenings.  Then, the discretized problem on $V_k$ is to find
$u_k \in V_k$ such that
\begin{equation}
\label{eqn:galerkinForm}
a(u_{k},v_k) = (f,v_k)_{L^2(\Omega)}, \ \ \ \forall v_k \in V_k.
\end{equation}
This problem has a unique solution, and the sequence
$\{u_k\}$ converges to $u^*$ \cite{BrennerScott94}.
%
The $L^2$-projection of the linear operator corresponding to the
bilinear form $a(\cdot\,,\cdot)$ onto $V_k$ is defined as the linear
operator $A_k : V_{k} \rightarrow V_{k}$ such that
\begin{equation}
\label{eqn:fematDef}
(A_{k} v_k,w_k)_{L^2(\Omega)} = a(v_k,w_k),  \ \ \ \forall v_k,w_k \in V_k.
\end{equation}
The operator $A_k$ is self-adjoint
with respect to the $L^2$-inner product and positive definite.
Let $\{\phi_1^k,\phi_2^k,\ldots,\phi_{N_k}^k\}$ be a basis for $V_k$ and
denote by $\mathbf{A_k}$ the representation of $A_k$ in that
basis. Then, \eqref{eqn:fematDef} becomes the linear matrix equation
for the coefficient vector $\mathbf{u}_k\in \mathbb{R}^{N_k}$
\begin{equation}
\mathbf{A}_k\mathbf{u}_k = \mathbf{f}_k,
\end{equation}
where, for $i,j = 1,2,\ldots,N_k$, the components of $\mathbf A_k$, $\mathbf u_k$ and $\mathbf f_k$ are given by
\begin{align*}
(\mathbf{A}_k)_{ij} =& a(\phi_i^k,\phi_j^k), \\
(\mathbf{f}_{k})_j   =& (f,\phi_j^k)_{L^2(\Omega)}, \\
(\mathbf{M}_k)_{ij} =& (\phi_i^k,\phi_j^k)_{L^2(\Omega)},
\end{align*}
where the integrals on the right hand sides are often approximated
using numerical quadrature.
Here, $\mathbf{M}_k$ is the mass matrix, which appears in the
approximation of the $L^2$-inner product in $V_k$ since
$(u_k,v_k)_{L^2(\Omega)} = \mathbf{u}_k^T\mathbf{M}_k\mathbf{v}_k$ for
all $u_k,v_k\in V_k$ with corresponding coefficient vectors
$\mathbf{u}_k,\mathbf{v}_k\in \mathbb{R}^{N_k}$.


\subsection{Restriction and prolongation} 
\label{sub:restriction_&_prolongation}
Since the coarse-grid space is a subspace of the fine-grid
space, any coarse-grid function $v_{k-1}$ can be expanded in
terms of the fine-grid basis functions,
\begin{equation} 
  v_{k-1} = \sum_{i=1}^{N_{k-1}} \mathbf v_{i,k-1}\phi_i^{k-1} = \sum_{j=1}^{N_k} \mathbf v_{j,k}\phi_j^k, 
\end{equation} 
where, $\mathbf v_{i,k}$ and $\mathbf v_{i,k-1}$ are the coefficients in the basis
expansion for $v_{k-1}$ on the fine and coarse grids, respectively.

The application of the prolongation operator can be represented as a matrix-vector
product with the input vector as the coarse grid nodal values and the
output as the fine grid nodal values \cite{SampathBiros10}. The matrix
entries of this operator are thus the coarse grid shape functions evaluated at the fine-grid
vertices, $p_i$, i.e.,
\begin{equation}
	\label{eq:Pstencil}
	\mathbf P_{\!ij} = \phi_j^{k-1}(p_i) \quad \text{ for } 1\le i \le N_k, 1\le j\le N_{k-1}. 
\end{equation}
This gives rise to two different operators depending on whether the
coarse grid is obtained via $h$-coarsening or whether it is obtained
via $p$-coarsening; see Figure~\ref{fig:approaches} for an
illustration of the two cases. The restriction operator is the adjoint
of the prolongation operator with respect to the mass-weighted inner
products. This only requires the application of the transpose of the
prolongation operator to vectors. 


\subsection{Meshing and implementation} 
\label{sub:meshing}

For the numerical comparisons in this work we consider domains that
are the image of a square or a cube under a
diffeomorphism, i.e., a smooth mapping from the reference domain
$S\coloneqq[0,1]^d$ to the physical domain $\Omega$. Hexahedral finite
element meshes and tensorized nodal basis function based on
Legende-Gauss-Lobotto (LGL) points are used.  We use isoparametric
elements to approximate the geometry of $\Omega$, i.e., on each element the geometry
diffeomorphism is approximated using the same basis
functions as the
finite element approximation. The Jacobians for this transformation
are computed at every quadrature point, and Gauss quadrature is used
to numerically approximate integrals.
We assume that the coefficient $\mu$ is a given function, which, at each level, can be evaluated at the respective quadrature points.
We restrict our comparisons to uniformly refined conforming meshes and our implementation, written in Matlab, is publicly available. It allows comparisons of different smoothing and
coarsening methods for high-order discretized problems in two and three dimensions, and can
easily be modified or extended. It does not support distributed memory
parallelism, and is
restricted to conforming meshes that can be mapped to a square (in 2D) or a
cube (in 3D). While, in practice, matrix assembly for
high-order discretizations is discouraged, we use sparse
assembled operators in this prototype implementation.


Note that for hexahedral elements in combination with a
tensorial finite element basis, the effect of matrix-free operations
for higher-order elements can be quite significant\footnote{For tetrahedral elements,
  this difference might be less pronounced.} in terms of
floating point operations, memory requirements, and actual run time:

\begin{itemize}
\item \emph{Memory requirements for assembled matrices:} For an order
  $p$, assembled element matrices are dense and of the size
  $(p+1)^3\times(p+1)^3$. For $p=9$, for instance, $(p+1)^3=1000$ and
  thus each element contributes $10^6$ entries to the assembled
  stiffness matrix, and each row in the matrix contains, on average,
  several 1000 nonzero entries. Thus, for high orders, memory becomes
  a significant issue.
\item \emph{Floating point operations for matrix-free versus assembled
  MatVec:} For hexahedral elements, the operation count for a
  tensorized matrix-free matvec is $\mathcal{O}(p^4)$ as opposed to
  $\mathcal{O}(p^6)$ for a fully assembled matrix \cite{orszag80,DevilleFischerMund02}.
\end{itemize}
  
  Detailed theoretical and experimental arguments in favor of matrix-free approaches, especially for high-order discretizations can be found in \cite{DevilleFischerMund02,BursteddeGhattasGurnisEtAl08,MayBrownLePourhiet14}





\section{Multigrid approaches for high-order finite element discretizations}
\label{sec:approaches}

In this section, we summarize different multigrid approaches
for high-order/spectral finite element
discretizations, which can either be used as
a solver or can serve as a preconditioner within a Krylov method. We
summarize different approaches for the construction of
multilevel hierarchies in~\S\ref{subsec:hierarchy} and discuss
smoothers in~\S\ref{subsec:smoothers}.

\subsection{Hierarchy construction, restriction and prolongation operators}\label{subsec:hierarchy}
There are several possibilities to build a multilevel hierarchy for
high-order discretized problems; see Figure~\ref{fig:approaches}. One
option is the construction of a geometric mesh hierarchy while keeping
the polynomial order unchanged; we refer to this approach as
high-order \emph{$h$-multigrid}. An alternative is to construct coarse
problems by reducing the polynomial degree of the finite element basis
functions, possibly followed by standard geometric multigrid; this is
commonly referred to as \emph{$p$-multigrid}. For unstructured
high-order element discretizations, where geometric coarsening is
challenging, using an algebraic multigrid hierarchy of a
\emph{low-order approximation to the high-order operator} as a
preconditioner has proven efficient. Some details of
these different approaches are summarized next.

\subsubsection{$h$-multigrid}\label{subsec:h}
A straightforward extension of low-order to high-order geometric
multigrid is to use the high-order discretization of the operator for
the residual computation on each multigrid level, combined with
high-order restriction and prolongation operators
(see \S\ref{sub:restriction_&_prolongation}).
%
For hexahedral (or quadrilateral) meshes, the required high-order residual
computations and the application of the interpolation and restriction operators
can often be accelerated using elementwise computations and
tensorized finite element basis functions, as is common
in spectral element methods \cite{DevilleFischerMund02}.

\subsubsection{$p$-multigrid}\label{subsec:p}
In the $p$-multigrid approach, a multigrid hierarchy is obtained by reducing
the polynomial order of the element basis functions.
Starting from an order-$p$ polynomial basis (for
simplicity, we assume here that $p$ is a power of 2), the coarser
grids correspond to polynomials of order $p/2, p/4,\ldots,1$, followed
by geometric coarsening of the $p=1$ grid (i.e., standard low-order
geometric multigrid).  As for high-order $h$-multigrid, devising
smoothers can be a challenge for $p$-multigrid.  Moreover, one often
finds dependence of the convergence factor on the order of the
polynomial basis \cite{MadayMunoz89}.

\subsubsection{Preconditioning by lower-order operator} \label{subsec:low}
In this defect correction approach (see
\cite{TrottenbergOosterleeSchuller01, Hackbusch85}), the high-order
residual is iteratively corrected using a low-order operator,
obtained by overlaying the high-order nodes with a low-order
(typically linear) finite element mesh. While the resulting low-order
operator has the same number of unknowns as the high-order operator,
it is much sparser and can, thus, be assembled efficiently and
provided as input to an algebraic multigrid method, which computes a
grid hierarchy through algebraic point aggregation.  This construction
of a low-order preconditioner based on the nodes of the high-order
discretization is used, for instance in~\cite{Brown10, Kim07,
  DevilleMund90, HeysManteuffelMcCormickEtAl05}. Due to the black-box
nature of algebraic multigrid, it is particularly attractive for
high-order discretizations on unstructured meshes. Note that even if
the mesh is structured, it is not
straightforward to use low-order geometric multigrid since the nodes---inherited from the high-order discretization---are not evenly spaced
(see Figure~\ref{fig:approaches}).

\subsection{Smoothers}\label{subsec:smoothers}
In our numerical comparisons, we focus on point smoothers but we also compare
with results obtained with an elementwise block-Jacobi smoother.
In this section, we summarize different
smoothers and numerically study their behavior for high-order
discretizations.  Note that, multigrid
smoothers must target the reduction of the error components in the
upper half of the spectrum.

\subsubsection{Point smoothers}
We compare the Jacobi and the symmetric successive over
relaxation (SSOR) smoothers, as well as a Chebyshev-accelerated Jacobi
smoother~\cite{Brandt77}. All of these smoothers require the diagonal of the
system matrix; if this matrix is not  assembled (i.e., in a matrix-free approach),
these diagonal entries must be computed in a setup step; for high-order discretizations on deformed meshes, this can be a significant computation.  Note that the
parallelization of Gauss-Seidel smoothers (such as SSOR) requires coloring of
unknowns in parallel, and, compared to Jacobi smoothing, more
complex communication in a distributed memory implementation. The
Chebyshev-accelerated Jacobi method is an alternative to SSOR; it
can significantly improve over Jacobi smoothing, while being as simple to
implement~\cite{AdamsBrezinaHuEtAl03}. The acceleration of Jacobi smoothing
with Chebyshev polynomials requires knowledge of the maximum eigenvalue of
the system matrix, usually estimated during setup with an iterative solver.


\subsubsection{Comparison of point smoothers}\label{subsubsec:ptsmoothcomparison}
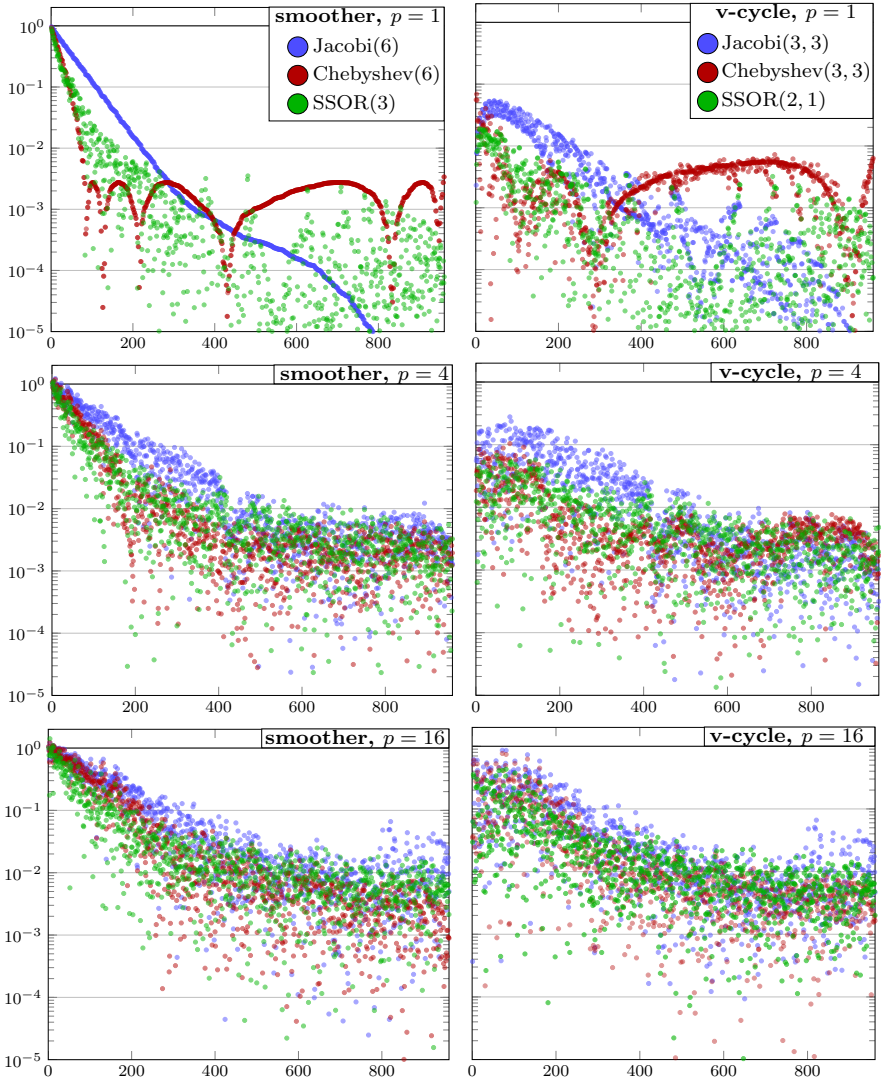
\begin{figure}
	\centering
\IfFileExists{homg-figure1.pdf}{\includegraphics[width=0.49\textwidth]{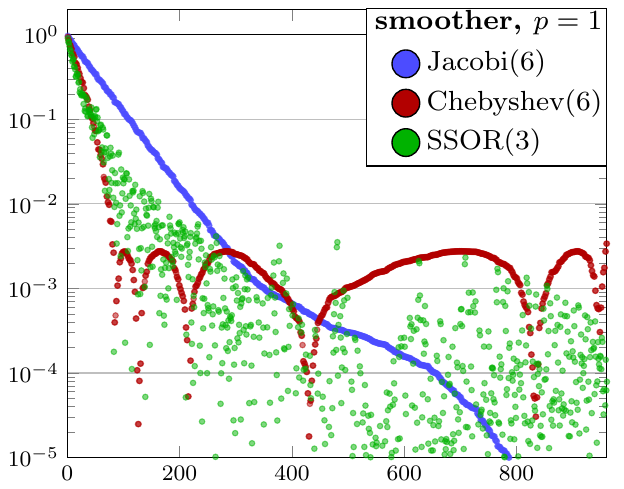}}{
		\begin{tikzpicture}[scale=0.9]
		\begin{semilogyaxis}[ymajorgrids,ymin=1e-5,ymax=2,xmin=0,xmax=961]
		\addplot[color=black]  table[x=dof, y=u]{smoother-const-box.dat};
		\addplot[color=blue!70, opacity=0.5,only marks, mark=*,mark size=1pt]   table[x=dof, y=jacobi1]{smoother-const-box.dat};
		\addplot[color=red!70!black, opacity=0.5,only marks, mark=*,mark size=1pt] table[x=dof, y=chebyshev1]{smoother-const-box.dat};
		\addplot[color=green!70!black,only marks, opacity=0.5,mark=*,mark size=1pt]  table[x=dof, y=ssor1]{smoother-const-box.dat};
		\end{semilogyaxis}
		\draw[black, fill=white] (3.8, 3.7) rectangle +(3.05,1.7);
		\node at (5.35, 5.53) {\bf \small{smoother, $p=1$}}; 
		\node[fill=blue!70, draw, circle,minimum width=0.1cm] at (4.3, 5.0) {}; 
		\node[fill=red!70!black, draw, circle,minimum width=0.1cm] at (4.3, 4.5) {};
		\node[fill=green!70!black, draw, circle,minimum width=0.1cm] at (4.3, 4.0) {};
		\node[text width=1.9cm] at (5.75, 5.0) {\small Jacobi$(6)$};
		\node[text width=1.9cm] at (5.75, 4.5) {\small Chebyshev$(6)$};
		\node[text width=1.9cm] at (5.75, 4.0) {\small SSOR$(3)$};
		\end{tikzpicture}
		}
		\IfFileExists{homg-figure2.pdf}{\includegraphics[width=0.45\textwidth]{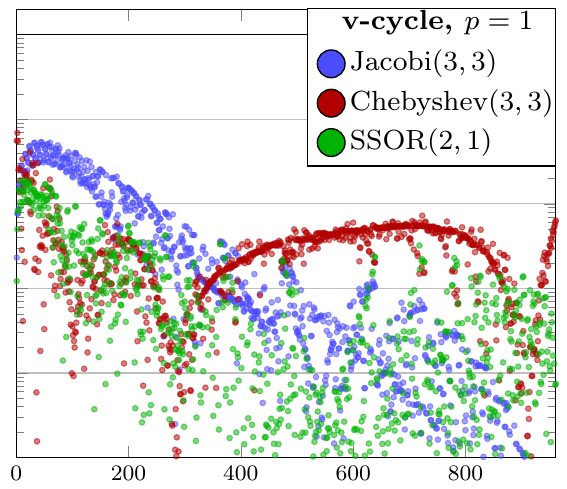}}{
		\begin{tikzpicture}[scale=0.9]
		\begin{semilogyaxis}[ymajorgrids,ymin=1e-5,ymax=2,xmin=0,xmax=961,yticklabels={,,}]
		\addplot[color=black]  table[x=dof, y=u]{vcycle-const-box.dat};
		\addplot[color=blue!70,opacity=0.5,only marks, mark=*,mark size=1pt]   table[x=dof, y=jacobi1]{vcycle-const-box.dat};
		\addplot[color=red!70!black,opacity=0.5,only marks, mark=*,mark size=1pt] table[x=dof, y=chebyshev1]{vcycle-const-box.dat};
		\addplot[color=green!70!black,opacity=0.5,only marks, mark=*,mark size=1pt]  table[x=dof, y=ssor1]{vcycle-const-box.dat};
		\end{semilogyaxis}
		\draw[black, fill=white] (3.7, 3.7) rectangle +(3.15,1.7);
		\node at (4.6, 5.53) {\bf \small{two-grid correction, $p=1$}}; 
		\node[fill=blue!70, draw, circle,minimum width=0.1cm] at (4.0, 5.0) {}; 
		\node[fill=red!70!black, draw, circle,minimum width=0.1cm] at (4.0, 4.5) {};
		\node[fill=green!70!black, draw, circle,minimum width=0.1cm] at (4.0, 4.0) {};
		\node[text width=2.1cm] at (5.55, 5.0) {\small Jacobi$(3,3)$};
		\node[text width=2.1cm] at (5.55, 4.5) {\small Chebyshev$(3,3)$};
		\node[text width=2.1cm] at (5.55, 4.0) {\small SSOR$(2,1)$};
		\end{tikzpicture}
		}
	\\
	\IfFileExists{homg-figure3.pdf}{\includegraphics[width=0.49\textwidth]{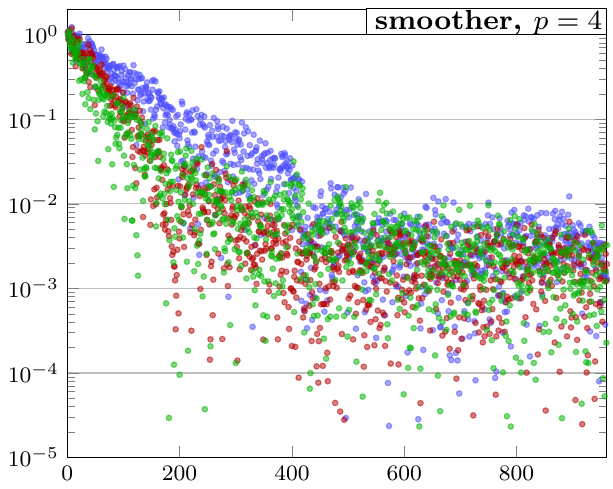}}{
		\begin{tikzpicture}[scale=0.9]
		\begin{semilogyaxis}[ymajorgrids,ymin=1e-5,ymax=2,xmin=0,xmax=961]
		\addplot[color=black]  table[x=dof, y=u]{smoother-const-box.dat};
		\addplot[color=blue!70,opacity=0.5,only marks, mark=*,mark size=1pt]   table[x=dof, y=jacobi4]{smoother-const-box.dat};
		\addplot[color=red!70!black,opacity=0.5,only marks, mark=*,mark size=1pt] table[x=dof, y=chebyshev4]{smoother-const-box.dat};
		\addplot[color=green!70!black,opacity=0.5,only marks, mark=*,mark size=1pt]  table[x=dof, y=ssor4]{smoother-const-box.dat};
		\end{semilogyaxis}
		\node at (5.35, 5.53) {\bf \small{smoother, $p=4$}};
		\end{tikzpicture}
		}
		\IfFileExists{homg-figure4.pdf}{\includegraphics[width=0.45\textwidth]{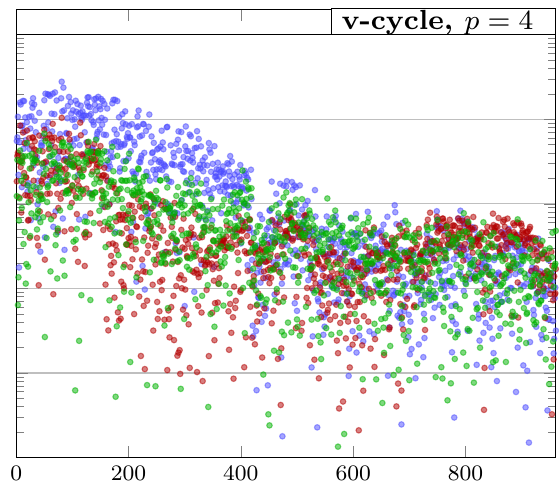}}{
		\begin{tikzpicture}[scale=0.9]
		\begin{semilogyaxis}[ymajorgrids,ymin=1e-5,ymax=2,xmin=0,xmax=961,yticklabels={,,}]
		\addplot[color=black]  table[x=dof, y=u]{vcycle-const-box.dat};
		\addplot[color=blue!70,only marks,opacity=0.5, mark=*,mark size=1pt]   table[x=dof, y=jacobi4]{vcycle-const-box.dat};
		\addplot[color=red!70!black,only marks,opacity=0.5, mark=*,mark size=1pt] table[x=dof, y=chebyshev4]{vcycle-const-box.dat};
		\addplot[color=green!70!black,only marks,opacity=0.5, mark=*,mark size=1pt]  table[x=dof, y=ssor4]{vcycle-const-box.dat};
		\end{semilogyaxis}
			\node at (4.6, 5.53) {\bf \small{two-grid correction, $p=4$}};
		\end{tikzpicture}
		}
	\\
	\IfFileExists{homg-figure5.pdf}{\includegraphics[width=0.49\textwidth]{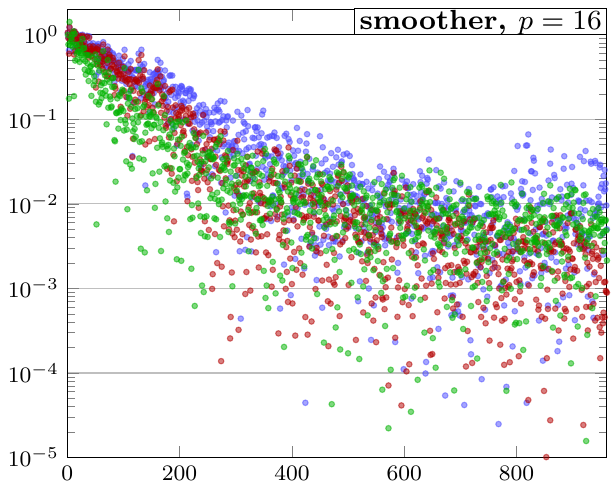}}{
		\begin{tikzpicture}[scale=0.9]
		\begin{semilogyaxis}[ymajorgrids,ymin=1e-5,ymax=2,xmin=0,xmax=961]
		\addplot[color=black]  table[x=dof, y=u]{smoother-const-box.dat};
		\addplot[color=blue!70,opacity=0.5,only marks, mark=*,mark size=1pt]   table[x=dof, y=jacobi16]{smoother-const-box.dat};
		\addplot[color=red!70!black,opacity=0.5,only marks, mark=*,mark size=1pt] table[x=dof, y=chebyshev16]{smoother-const-box.dat};
		\addplot[color=green!70!black,opacity=0.5,only marks, mark=*,mark size=1pt]  table[x=dof, y=ssor16]{smoother-const-box.dat};
		\end{semilogyaxis}
		\node at (5.25, 5.53) {\bf \small{smoother, $p=16$}};
		\end{tikzpicture}
		}
		\IfFileExists{homg-figure6.pdf}{\includegraphics[width=0.45\textwidth]{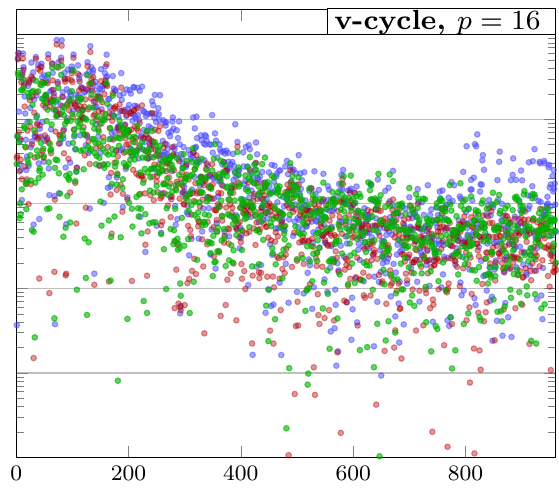}}{
		\begin{tikzpicture}[scale=0.9]
		\begin{semilogyaxis}[ymajorgrids,ymin=1e-5,ymax=2,xmin=0,xmax=961,yticklabels={,,}]
		\addplot[color=black]  table[x=dof, y=u]{vcycle-const-box.dat};
		\addplot[color=blue!70,opacity=0.5,only marks, mark=*,mark size=1pt]   table[x=dof, y=jacobi16]{vcycle-const-box.dat};
		\addplot[color=red!70!black,opacity=0.4,only marks, mark=*,mark size=1pt] table[x=dof, y=chebyshev16]{vcycle-const-box.dat};
		\addplot[color=green!70!black,opacity=0.6,only marks, mark=*,mark size=1pt]  table[x=dof, y=ssor16]{vcycle-const-box.dat};
		\end{semilogyaxis}
		\node at (4.6, 5.53) {\bf \small{two-grid correction, $p=16$}};
		\end{tikzpicture}
		}
	\caption{\label{fig:smoothers} Error decay for different point
          smoothers when used as solver (left column) and when used in
          a single two-grid step with exact coarse grid solution
          (right column) for a two-dimensional, constant coefficient
          Laplace problem on a unit square (problem {\bf 2d-const}
          specified in \S\ref{subsec:tests}). To keep the
          number of unknowns the same accross all polynomial orders,
          meshes of $32\times 32$, $8\times 8$ and $2\times 2$
          elements are used for polynomial orders $p=1$, $p=4$ and
          $p=16$, respectively.  The horizontal axis is the
          index for the eigenvectors of the system matrix $\mathbf{A}_k$, and
          the vertical axis is the magnitude of the error component
          for each eigenvector. The eigenvectors are ordered such that
          the corresponding eigenvalues are ascending; thus, due to
          the properties of $\mathbf{A}_k$, the
          smoothness in every eigenvector decays from left to
          right.  The system right hand side is zero and the
          initialization is chosen to have all unit coefficients in
          the eigenvector expansion. A total of six smoothing steps is
          used for all methods, and the coarse problem in the two-grid
          step is solved by a direct solver.}
\end{figure}

\begin{figure}
	\centering
	\IfFileExists{homg-figure7.pdf}{\includegraphics[width=0.49\textwidth]{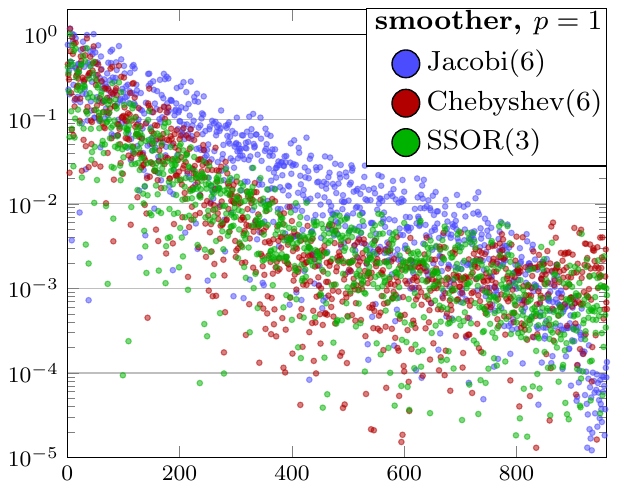}}{
		\begin{tikzpicture}[scale=0.9]
		\begin{semilogyaxis}[ymajorgrids,ymin=1e-5,ymax=2,xmin=0,xmax=961]
		\addplot[color=black]  table[x=dof, y=u]{smoother-var-shell.dat};
		\addplot[color=blue!70, opacity=0.5,only marks, mark=*,mark size=1pt]   table[x=dof, y=jacobi1]{smoother-var-shell.dat};
		\addplot[color=red!70!black, opacity=0.5,only marks, mark=*,mark size=1pt] table[x=dof, y=chebyshev1]{smoother-var-shell.dat};
		\addplot[color=green!70!black,only marks, opacity=0.5,mark=*,mark size=1pt]  table[x=dof, y=ssor1]{smoother-var-shell.dat};
		\end{semilogyaxis}
		\draw[black, fill=white] (3.8, 3.7) rectangle +(3.05,1.65);
		\node at (5.35, 5.53) {\bf \small{smoother, $p=1$}}; 
		\node[fill=blue!70, draw, circle,minimum width=0.1cm] at (4.3, 5.0) {}; 
		\node[fill=red!70!black, draw, circle,minimum width=0.1cm] at (4.3, 4.5) {};
		\node[fill=green!70!black, draw, circle,minimum width=0.1cm] at (4.3, 4.0) {};
		\node[text width=1.9cm] at (5.75, 5.0) {\small Jacobi$(6)$};
		\node[text width=1.9cm] at (5.75, 4.5) {\small Chebyshev$(6)$};
		\node[text width=1.9cm] at (5.75, 4.0) {\small SSOR$(3)$};
		\end{tikzpicture}
		}
		\IfFileExists{homg-figure8.pdf}{\includegraphics[width=0.45\textwidth]{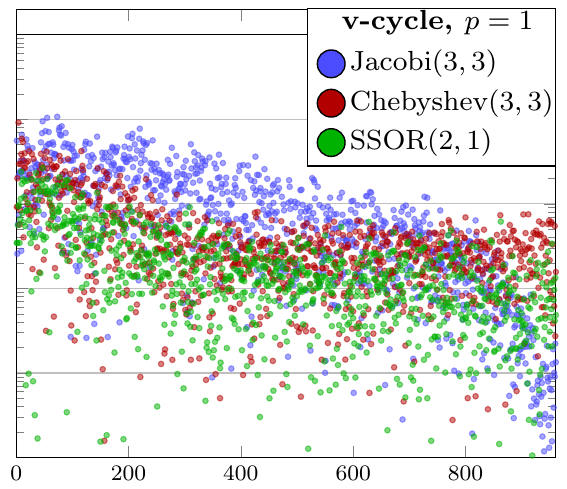}}{
		\begin{tikzpicture}[scale=0.9]
		\begin{semilogyaxis}[ymajorgrids,ymin=1e-5,ymax=2,xmin=0,xmax=961,yticklabels={,,}]
		\addplot[color=black]  table[x=dof, y=u]{vcycle-var-shell.dat};
		\addplot[color=blue!70,opacity=0.5,only marks, mark=*,mark size=1pt]   table[x=dof, y=jacobi1]{vcycle-var-shell.dat};
		\addplot[color=red!70!black,opacity=0.5,only marks, mark=*,mark size=1pt] table[x=dof, y=chebyshev1]{vcycle-var-shell.dat};
		\addplot[color=green!70!black,opacity=0.5,only marks, mark=*,mark size=1pt]  table[x=dof, y=ssor1]{vcycle-var-shell.dat};
		\end{semilogyaxis}
		\draw[black, fill=white] (3.7, 3.7) rectangle +(3.15,1.65);
		\node at (4.65, 5.53) {\bf \small{two-grid correction, $p=1$}}; 
		\node[fill=blue!70, draw, circle,minimum width=0.1cm] at (4.0, 5.0) {}; 
		\node[fill=red!70!black, draw, circle,minimum width=0.1cm] at (4.0, 4.5) {};
		\node[fill=green!70!black, draw, circle,minimum width=0.1cm] at (4.0, 4.0) {};
		\node[text width=2.1cm] at (5.55, 5.0) {\small Jacobi$(3,3)$};
		\node[text width=2.1cm] at (5.55, 4.5) {\small Chebyshev$(3,3)$};
		\node[text width=2.1cm] at (5.55, 4.0) {\small SSOR$(2,1)$};
		\end{tikzpicture}
		}
	\\
	\IfFileExists{homg-figure9.pdf}{\includegraphics[width=0.49\textwidth]{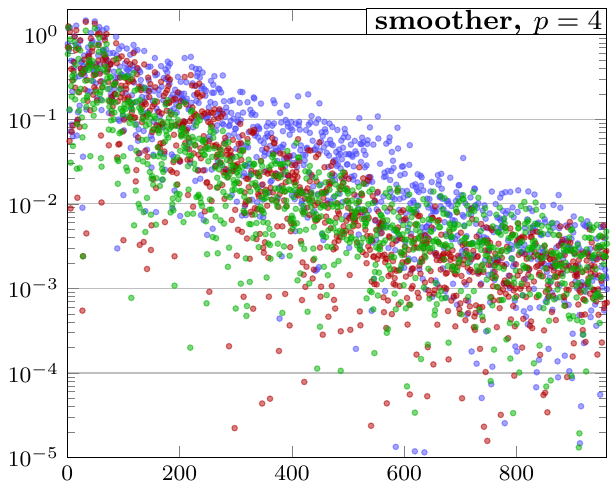}}{
		\begin{tikzpicture}[scale=0.9]
		\begin{semilogyaxis}[ymajorgrids,ymin=1e-5,ymax=2,xmin=0,xmax=961]
		\addplot[color=black]  table[x=dof, y=u]{smoother-var-shell.dat};
		\addplot[color=blue!70,opacity=0.5,only marks, mark=*,mark size=1pt]   table[x=dof, y=jacobi4]{smoother-var-shell.dat};
		\addplot[color=red!70!black,opacity=0.5,only marks, mark=*,mark size=1pt] table[x=dof, y=chebyshev4]{smoother-var-shell.dat};
		\addplot[color=green!70!black,opacity=0.5,only marks, mark=*,mark size=1pt]  table[x=dof, y=ssor4]{smoother-var-shell.dat};
		\end{semilogyaxis}
		\node at (5.35, 5.53) {\bf \small{smoother, $p=4$}};
		\end{tikzpicture}
		}
		\IfFileExists{homg-figure10.pdf}{\includegraphics[width=0.45\textwidth]{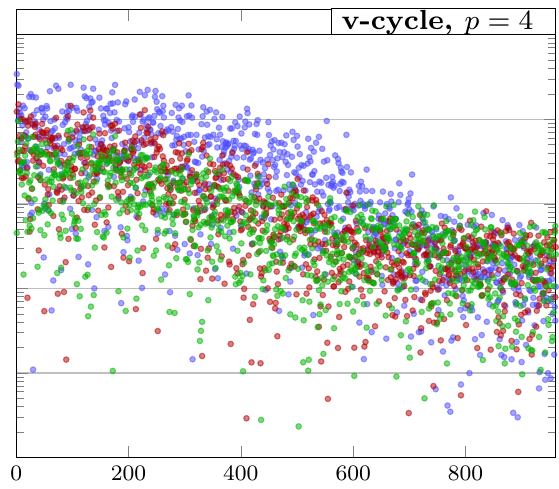}}{
		\begin{tikzpicture}[scale=0.9]
		\begin{semilogyaxis}[ymajorgrids,ymin=1e-5,ymax=2,xmin=0,xmax=961,yticklabels={,,}]
		\addplot[color=black]  table[x=dof, y=u]{vcycle-var-shell.dat};
		\addplot[color=blue!70,only marks, mark=*,opacity=0.5,mark size=1pt]   table[x=dof, y=jacobi4]{vcycle-var-shell.dat};
		\addplot[color=red!70!black,only marks, mark=*,opacity=0.5,mark size=1pt] table[x=dof, y=chebyshev4]{vcycle-var-shell.dat};
		\addplot[color=green!70!black,only marks, mark=*,opacity=0.5,mark size=1pt]  table[x=dof, y=ssor4]{vcycle-var-shell.dat};
		\end{semilogyaxis}
		\node at (4.65, 5.53) {\bf \small{two-grid correction, $p=4$}};
		\end{tikzpicture}
		}
	\\
	\IfFileExists{homg-figure11.pdf}{\includegraphics[width=0.49\textwidth]{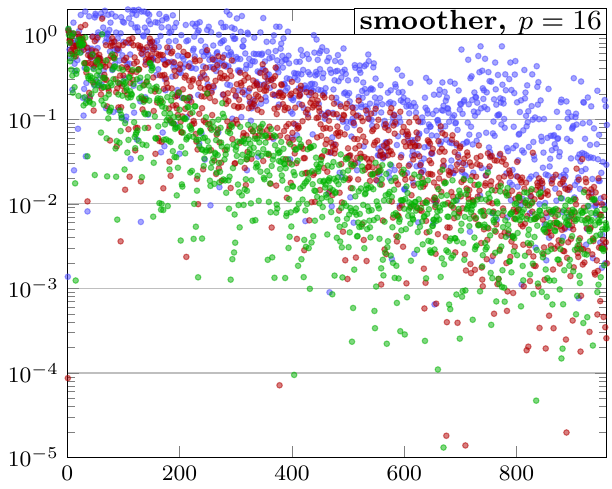}}{
		\begin{tikzpicture}[scale=0.9]
		\begin{semilogyaxis}[ymajorgrids,ymin=1e-5,ymax=2,xmin=0,xmax=961]
		\addplot[color=black]  table[x=dof, y=u]{smoother-var-shell.dat};
		\addplot[color=blue!70,opacity=0.5,only marks, mark=*,mark size=1pt]   table[x=dof, y=jacobi16]{smoother-var-shell.dat};
		\addplot[color=red!70!black,opacity=0.5,only marks, mark=*,mark size=1pt] table[x=dof, y=chebyshev16]{smoother-var-shell.dat};
		\addplot[color=green!70!black,opacity=0.5,only marks, mark=*,mark size=1pt]  table[x=dof, y=ssor16]{smoother-var-shell.dat};
		\end{semilogyaxis}
		\node at (5.25, 5.53) {\bf \small{smoother, $p=16$}};
		\end{tikzpicture}
		}
		\IfFileExists{homg-figure12.pdf}{\includegraphics[width=0.45\textwidth]{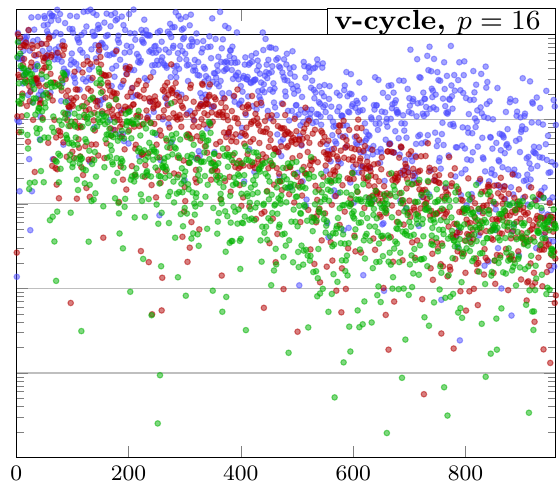}}{
		\begin{tikzpicture}[scale=0.9]
		\begin{semilogyaxis}[ymajorgrids,ymin=1e-5,ymax=2,xmin=0,xmax=961,yticklabels={,,}]
		\addplot[color=black]  table[x=dof, y=u]{vcycle-var-shell.dat};
		\addplot[color=blue!70,opacity=0.5,only marks, mark=*,mark size=1pt]   table[x=dof, y=jacobi16]{vcycle-var-shell.dat};
		\addplot[color=red!70!black,opacity=0.5,only marks, mark=*,mark size=1pt] table[x=dof, y=chebyshev16]{vcycle-var-shell.dat};
		\addplot[color=green!70!black,opacity=0.5,only marks, mark=*,mark size=1pt]  table[x=dof, y=ssor16]{vcycle-var-shell.dat};
		\end{semilogyaxis}
		\node at (4.65, 5.53) {\bf \small{two-grid correction, $p=16$}};
		\end{tikzpicture}
		}
	\caption{\label{fig:smoothers-var} Shown is the same comparison as
          in Figure~\ref{fig:smoothers}, but for the two-dimensional
          warped geometry, variable coefficient problem {\bf 2d-var}
          specified in \S\ref{subsec:tests}.}
\end{figure}

In Figures~\ref{fig:smoothers} and~\ref{fig:smoothers-var}, we compare
the efficiency of these point smoothers for different polynomial
orders and constant and varying coefficients. For that purpose, we
compute the eigenvectors of the system matrix, choose a
zero right hand side and an initialization that has all unit
coefficients in the basis given by these eigenvectors. For the
polynomial orders $p=1,4,16$, we compare the performance of point
smoothers with and without a 2-level v-cycle with exact coarse solve.
The coarse grid for all polynomial orders is obtained using $h$-coarsening.
We depict the
coefficients after six smoothing steps in the left column, and the results
obtained for a two-grid method\footnote{For simplicity, we chose two
  grids in our tests; the results for a multigrid v-cycle are
  similar.} with three pre- and three post-smoothing steps (and thus overall six
smoothing steps on the finest grid) in the right column. The SSOR smoother
uses a lexicographic ordering of the unknowns, and we employ two pre- and one
post-smoothing steps, which again amounts to overall six smoothing steps
on the finest grid. The damping factors for Jacobi and SSOR smoothing
are $\omega = 2/3$ and $\omega=1$, respectively. The Chebyshev
smoother targets the part of the
spectrum given by $[\lambda_\text{max}/4,\lambda_\text{max}]$, where
$\lambda_\text{max}$ is the maximum eigenvalue of the system matrix, which is estimated using 10 iterations of the Arnoldi algorithm.

The results for the constant coefficient Laplacian operator on the
unit square (see Figure~\ref{fig:smoothers}) show that all point smoothers
decrease the error components in the upper half of the spectrum;
however,
the decrease is smaller for high-order elements. Observe that compared to Jacobi smoothing, Chebyshev accelerated Jacobi smoothing dampens a larger
part of the spectrum.  Both, the Chebyshev and
SSOR methods outperform Jacobi smoothing, in particular for higher
orders. Combining the smoothers with a two-grid cycle, all error
components are decreased for all smoothers (and thus the resulting
two-grid methods converge, see Table~\ref{tab:box} in \S\ref{subsec:results}), but the error decreases slower for
higher polynomial orders. For high polynomial orders, a two-grid
iteration with SSOR smoothing results in a much better error
reduction than Jacobi or Chebyshev smoothing.

In Figure~\ref{fig:smoothers-var}, we study the performance of
different smoothers for the test problem {\bf 2d-var}, defined in
\S\ref{subsec:tests}. In this problem, we solve
\eqref{eq:Poisson} with a strongly (but smoothly) varying coefficient
$\mu$ on a deformed domain $\Omega$.  Compared to the constant
coefficient case, Jacobi smoothing performs worse, both,
when used as a solver and as a smoother. Let us focus on the two-grid
correction for polynomial order $p=16$ and compare with the results
obtained when using multigrid as a solver, shown in Table~\ref{tab:2d-fan}.
Jacobi smoothing does not lead to a converging two-grid algorithm, as
several coefficients are amplified by the two-grid cycle. For
Chebyshev smoothing, the multigrid v-cycle converges slowly
although one or two coefficients appear amplified in the two-grid
iteration. This convergence can be explained by the fact that errors
can be interchanged between different eigenvectors in the v-cycle.
SSOR smoothing combined with the two-grid method retains a significant
error reduction rate and, as a consequence, converges quickly.

\subsubsection{Block-Jacobi smoothing}\label{subsubsec:schwarz}
An alternative smoothing approach for high-order discretizations is
based on local block solves.  Since for high polynomial orders many
unknowns lie in the element interiors, Schwarz-type domain
decomposition smoothers are promising. For instance, they are more
stable for anisotropic meshes than point smoothers. A main challenge
of Schwarz-type smoothers is that they require the solution of dense
local systems.  This is either done by using direct methods or
approximations that allow for a fast iterative solution on hexahedral
meshes \cite{LottesFischer05, FischerLottes05}.
In
\S\ref{sec:numerics}, we compare the performance of point
smoothers with an elementwise block Jacobi smoothing.

\section{Numerical results}\label{sec:numerics}
In this section, we present a comprehensive comparison of our
algorithms for the solution of high-order discretizations of
\eqref{eq:Poisson}.  After introducing our test problems in
\S\ref{subsec:tests}, we present a simple model for the computational
cost of the different approaches in terms of matrix-vector
applications in \S\ref{subsec:complexity}. In \S\ref{subsec:measures},
we specify settings and metrics for our comparisons. The results of
these comparisons are presented and discussed in
\S\ref{subsec:results}.

\subsection{Test problems}\label{subsec:tests}
\begin{figure}
	\includegraphics[width=0.48\textwidth]{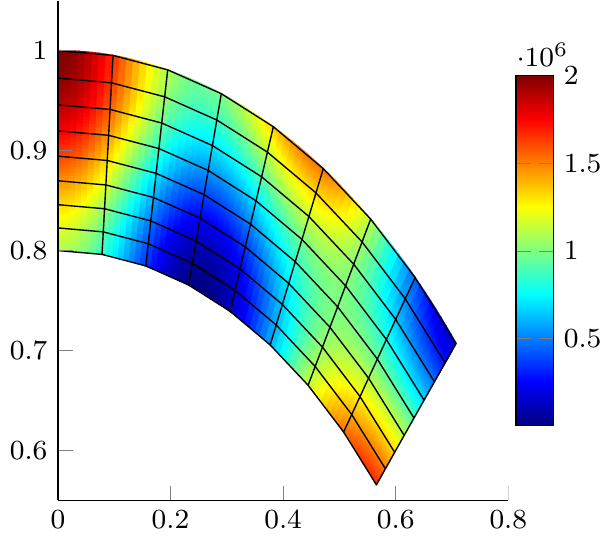}
	\includegraphics[width=0.48\textwidth]{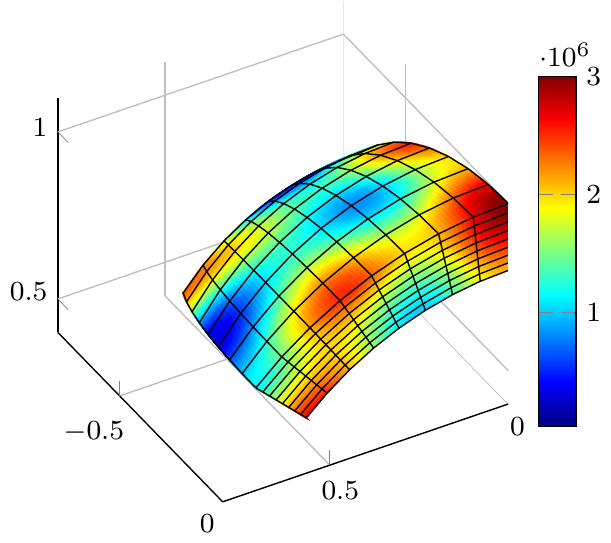}
	\caption{\label{fig:mesh} Two and three-dimensional warped
          meshes used in our numerical experiments. The color
          illustrates the logarithm of the coefficient field, which
          varies over six orders of magnitude.}
\end{figure}

We compare our algorithms for the solution of~\eqref{eq:Poisson} with
constant coefficient $\mu\equiv 1$ on the unit square and the unit cube,
and, with varying coefficients $\mu(\bs x)$, on the warped two and
three-dimensional domains shown in Figure~\ref{fig:mesh}. To be
precise, we consider the following four problems:
\begin{itemize}
\item {\bf 2d-const:} The domain $\Omega$ for the problem is the unit square, and $\mu\equiv
  1$.
\item {\bf 2d-var:} The warped two-dimensional domain $\Omega$ is
  shown on the left in Figure~\ref{fig:mesh}, and the varying
  coefficient is $\mu(x,y) = 1 + 10^6(\cos^2(2\pi x) + \cos^2(2\pi
  y))$. We also study a modification of this problem with a more
  oscillatory coefficient
  $\mu(x,y) = 1 + 10^6(\cos^2(10\pi x) +
  \cos^2(10\pi y))$, which we refer to as {\bf 2d-var$'$}.
\item {\bf 3d-const:} For this problem, $\Omega$ is the unit cube,
  and we use the constant coefficient $\mu\equiv 1$.
\item {\bf 3d-var:} The warped three-dimensional domain $\Omega$
  shown on the right of Figure~\ref{fig:mesh} is used; the
  varying coefficient is $\mu(x,y,z) = 1 + 10^6(\cos^2(2\pi x) +
  \cos^2(2\pi y) + \cos^2(2\pi z))$.
\end{itemize}

\subsection{Comparing the computational cost}\label{subsec:complexity}
To compare the computational cost of the different methods, we focus
on the matrix-vector multiplications on the finest multigrid level,
which dominate the overall computation. Denoting the number of
unknowns on the finest level by $N$, the computational cost---measured
in floating point operations (flops)---for a matrix-vector product is
$Ng_p$, where $g_p$ is the number of flops per unknown and the
subscript $p$ indicates the polynomial order used in the FEM basis.
Since high-order discretizations result in less sparse operators,
$g_1\le g_2\le \ldots$ holds. The actual value of $g_p$ depends
strongly on the implementation. Also note that the conversion from
$g_p$ to wall-clock time is not trivial, as wall-clock timings depend
on caching, vectorization, blocking and other effects. Thus,
although $g_p$ increases with $p$, wall-clock
times might not increase as significantly.
In general, high-order implementations allow more memory locality,
which often results in higher performance compared to low-order
methods.
This discussion, however, is beyond the scope of this paper.


The dominant computational cost per iteration of the high-order
multigrid approaches discussed in \S\ref{sec:approaches} can
thus be summarized as
\begin{equation}\label{eq:compcost}
  Ng_p(1+m(s_\text{pre}+s_\text{post})).
\end{equation}
Here, we denote by $s_\text{pre}$ and $s_\text{post}$ the number of
pre- and post-smoothing steps on the finest multigrid level,
respectively. Moreover, $m$ denotes the number of residual
computations (and thus matrix-vector computations) per smoothing step.
Jacobi smoothing and Chebyshev-accelerated Jacobi require $m=1$
matrix-vector multiplication per smoothing step, while SSOR requires
$m=2$ matrix-vector operations. If, in the approach discussed in
\S\ref{subsec:low}, the sparsified linear-element residual is used in the
smoother on the finest grid, the cost \eqref{eq:compcost} reduces to
\begin{equation}\label{eq:compcost2}
  N(g_p + g_1 m(s_\text{pre}+s_\text{post})).
\end{equation}
However, since the overall number of iterations increases (see
\S\ref{subsec:results}), this does not necessarily decrease the
solution time.


%
%

If the overall number of unknowns $N$ is kept fixed and the
solution is smooth, it is well known that the accuracy increases for
high-order discretizations. Due to the decreased sparsity of the
discretized operators, this does not automatically translate to more
accuracy per computation time; see, e.g.,~\cite{Brown10}. However, note
that many computations in, for instance, a multigrid
preconditioned conjugate gradient algorithm are of complexity
$\mathcal{O}(N)$ (see Algorithm~\ref{alg:pcg}) and are thus independent of
$g_p$. Thus, the computational cost of these steps does not depend on
the order of the discretization. Even if these $\mathcal{O}(N)$ steps
do not dominate the computation, they contribute to making high-order
discretizations favorable not only in terms of accuracy per unknown,
but also in terms of accuracy per computation time.

\begin{algorithm}[ht] 
  \caption{Complexity of individual steps in multigrid-preconditioned CG} \label{alg:pcg} 
  \begin{algorithmic}[1]
    \Require rhs and guess
    \Ensure  solution
    \While {not converged} 
    \State $\bs{h} = A \bs{p}$ 											\Comment $~~\quad\quad\quad\quad\mathcal{O}(Ng_p)$
    \State $\rho_r = (\rho, \bs{r})$								\Comment $~~\quad\quad\quad\quad\mathcal{O}(N)~~~$
    \State $\alpha = \rho_r / ( \bs{p}, \bs{h} )$		\Comment $~~\quad\quad\quad\quad\mathcal{O}(N)~~~$
    \State $\bs{u} = \bs{u} + \alpha\bs{p}$					\Comment $~~\quad\quad\quad\quad\mathcal{O}(N)~~~$
    \State $\bs{r} = \bs{r} - \alpha\bs{h}$					\Comment $~~\quad\quad\quad\quad\mathcal{O}(N)~~~$
    \State Convergence Test
    \State $\rho = M\bs{r}$ 												\Comment v-cycle $\quad\mathcal{O}(Ng_p)$
    \State $\beta = (\rho, \bs{r}) / \rho_r$				\Comment $~~\quad\quad\quad\quad\mathcal{O}(N)~~~$
    \State $\bs{p} = \rho + \beta\bs{p}$						\Comment $~~\quad\quad\quad\quad\mathcal{O}(N)~~~$
    \EndWhile
  \end{algorithmic}
\end{algorithm}

\subsection{Setup of comparisons}\label{subsec:measures}
We test the different multigrid schemes in two contexts: as solvers
and as preconditioners in a conjugate gradient (CG) method.  In tables
\ref{tab:box}--\ref{tab:3d-fan}, 
we report the number of multigrid v-cycles%
\footnote{each CG iteration uses a single multigrid
v-cycle as preconditioner} 
required to reduce the norm of the discrete
residual by a factor of $10^8$, where a ``-'' indicates that the
method did not converge within the specified maximum number of
iterations.
In particular, these tables report the following
information:
\begin{itemize}
\item[$\bullet$] The first column gives the polynomial \emph{order}
  used in the finite element discretization.
\item[$\bullet$] The columns labeled \emph{MG as solver} report
  the number of v-cycles required for convergence when multigrid is used as
  solver. The subcolums are:
  \begin{itemize}
  \item \emph{Jacobi(3,3)} denotes that 3 pre-smoothing and 3
    post-smoothing steps of a pointwise Jacobi smoother are used on
    each level. We use a damping factor $\omega=2/3$ in all experiments.
  \item \emph{Cheb(3,3)} indicates that Chebyshev-accelerated Jacobi
    smoothing is used, again using 3 pre-smoothing and 3
    post-smoothing steps. An estimate for the maximal eigenvalue of
    the linear systems on each level, as required by the Chebyshev
    method, is computed in a setup step using 10 Arnoldi iterations.
  \item \emph{SSOR(2,1)} denotes that a symmetric successive
    over-relaxation method is employed, with 2 pre-smoothing and 1
    post-smoothing steps. Note that each SSOR iteration amounts to a
    forward and a backward Gauss-Seidel smoothing step, and thus
    requires roughly double the computational work compared to Jacobi
    smoothing.  The SSOR smoother is based on a lexicographic ordering
    of the unknowns, and the damping factor is $\omega=1$.
  \end{itemize}
    For the two-dimensional problems reported in
    Tables~\ref{tab:box}--\ref{tab:2d-fan2}, we use a multigrid
    hierarchy with three levels corresponding to meshes with
    $32\times32$, $16\times16$ and $8\times8$ elements.
    The multigrid hierarchy for the three-dimensional
    tests reported in Tables~\ref{tab:3d-box} and \ref{tab:3d-fan}
    also has three levels with
    $8\times8\times8$, $4\times4\times4$ and $2\times2\times2$
    elements.  Note that for each smoother we report results
    for $h$-multigrid (columns marked by \emph{h}; see
    \S\ref{subsec:h}) as well as for $p$-multigrid (columns marked by
    \emph{p}; see \S\ref{subsec:p}). For $p$-multigrid, we restrict
    ourselves to orders that are powers of 2. After coarsening in $p$
    till $p=1$, we coarsen in $h$. For example, for the
    two-dimensional problems and $p=16$, we use a total of 7 grids;
    the first five all use meshes with $32\times32$ elements,
    and $p=16,8,4,2,1$, respectively, followed by two additional
    coarse grids of size $16\times16$ and $8\times8$, and
    $p=1$.
\item[$\bullet$] The columns labeled \emph{MG with pCG} present the
  number of conjugate gradient iterations required for the solution,
  where each iteration uses one multigrid v-cycle as preconditioner.
  The sub-columns correspond to different smoothers, as described
  above.
\item[$\bullet$] The columns labeled \emph{low-order MG pCG} report
  the number of CG iterations needed to solve the high-order system,
  when preconditioned with the low-order operator based on the
  high-order nodal points (see \S\ref{subsec:low}).  While in practice
  one would use algebraic multigrid to solve the linearized system
  approximately, in our tests we use a factorization method to solve
  the low-order system directly.  As a consequence, the reported
  iteration counts are a lower bound for the iteration counts one
  would obtain if the low-order system was inverted approximately by
  algebraic multigrid.

\end{itemize}

Note that the number of smoothing steps in the different methods
is chosen such that, for fixed polynomial order, the computational
work is comparable. Each multigrid v-cycle requires one residual
computation and overall six matrix-vector multiplications.  Following
the simple complexity estimates \eqref{eq:compcost} and
\eqref{eq:compcost2}, this amounts to a per-iteration cost of $7Ng_p$
for $h$- and $p$-multigrid, and of $N(g_1+6g_p)$ for the low-order
multigrid preconditioner. As a consequence, the iteration numbers
reported in the next section can be used to compare
the efficiency of the different methods.
Note that in our tests, we change the polynomial degree of the
finite element functions but retain the same mesh. This results in an
increasing number of unknowns as $p$ increases. Since, as illustrated
in \S\ref{subsec:num_mesh}, we observe mesh independent convergence
for fixed $p$, this does not influence the comparison.

\subsection{Summary of numerical results}\label{subsec:results}
Next, in \S\ref{subsec:num_point}, we compare the performance of
different point smoothers for the test problems presented in
\S\ref{subsec:tests}. Then, in \S\ref{subsec:num_mesh}, we
illustrate that the number of iterations is independent of the mesh
resolution. Finally, in \S\ref{subsec:num_block}, we study the
performance of a block Jacobi smoother for discretizations with
polynomial orders $p=8$ and $p=16$.

\begin{table}
  \caption{\label{tab:box} Iteration counts for the two-dimensional unit square
    problem {\bf 2d-const} defined in \S\ref{subsec:tests}. The finest
    mesh has $32\times 32$ elements and the multigrid hierarchy
    consists of  three meshes.  For $p$-multigrid, the
    polynomial order is first reduced to $p=1$, followed by two
    geometric coarsenings of the mesh.  For a detailed description of
    the different experiments reported in this table we refer to
    \S\ref{subsec:measures}.}  \centering
  \begin{tabular}{|r|c c|c c|c c||c c|c c|c c||c|} 
    \hline
    & \multicolumn{6}{c||}{MG as solver} & \multicolumn{6}{c||}{MG
      with pCG} & \!\!low-order MG\!\! \\
    \cline{2-13}
    \!\!\! order \!\!\!\! &  \multicolumn{2}{c|}{\!\!\scriptsize  Jacobi(3,3)\!\!} &  \multicolumn{2}{c|}{\!\!\scriptsize Cheb(3,3)\!\!} & \multicolumn{2}{c||}{\!\!\scriptsize  SSOR(2,1)\!\!} & \multicolumn{2}{c|}{\!\!\scriptsize Jacobi(3,3)\!\!} &  \multicolumn{2}{c|}{\!\!\scriptsize Cheb(3,3)\!\!} & \multicolumn{2}{c||}{\!\!\scriptsize SSOR(2,1)\!\!} & pCG \\
\hline
 & $h$ & $p$ & $h$ & $p$& $h$ & $p$& $h$ & $p$& $h$ & $p$& $h$ & $p$&
~ \\
 \cline{2-13}
1 & 6 & & 5 & & 5 & & 5 & & 4 & & 4 & & -   \\
2 & 7 & 7 & 5 & 6 & 5 & 5 & 5 & 5 & 4 & 4 & 4 & 4 & 14  \\
3 & 8 & & 6 & & 5 & & 6 & & 5 & & 4 & & 16   \\
4 & 9 & 8 & 6 & 6 & 5 & 5 & 6 & 5 & 5 & 5 & 4 & 4 & 16  \\
5 & 12 & & 8 & & 7 & & 7 & & 6 & & 5 & & 17 \\
6 & 12 & & 9 & & 7 & & 7 & & 6 & & 5 & & 18 \\
7 & 16 & & 12 & & 8 & & 8 & & 7 & & 6 & & 18  \\
8 & 17 & 14 & 13 & 10 & 8 & 7 & 9 & 8 & 7 & 6 & 6 & 5 & 19 \\
16 & 40 & 33 & 33 & 27 & 17 & 14 & 14 & 12 & 12 & 11 & 9 & 8 & 21 \\
\hline
  \end{tabular}
\end{table}

\subsubsection{Comparison of different multigrid/smoothing combinations}\label{subsec:num_point}
Tables \ref{tab:box}--\ref{tab:2d-fan2} present the number of
iterations obtained for various point smoothers and different
polynomial orders for the two-dimensional test problems. As can be
seen in Table~\ref{tab:box}, for {\bf 2d-const}
all solver variants converge in a
relatively small number of iterations for all polynomial
orders. However, the number of
iterations increases with the polynomial order $p$, in particular when
multigrid is used as a solver. Using multigrid as a preconditioner in
the conjugate gradient method results in a reduction of overall
multigrid v-cycles, in some cases even by a factor or two. Also, we
observe that SSOR smoothing generally performs better than the two
Jacobi-based smoothers. We find that the linear-order operator based
on the high-order nodes is a good preconditioner for the high-order
system. Note that if algebraic multigrid is used for the solution of
the low-order approximation, the smoother on the finest level can
either use the residual of the low-order or of the high-order
operator. Initial tests that mimic the use of the high-order residual
in the fine-grid smoother show that this has the potential to reduce the
number of iterations.

%

Let us now contrast these observations with the results for the
variable coefficient problems {\bf 2d-var} and {\bf 2d-var$'$}
summarized in Tables~\ref{tab:2d-fan} and \ref{tab:2d-fan2}. First,
note that all variants of the solver perform reasonably for
discretizations up to order $p=4$. When used as a solver, multigrid
either diverges or converges slowly for orders $p>4$. Convergence
is reestablished when multigrid is combined with CG. Using multigrid
with SSOR smoothing as preconditioner in CG yields, for orders up to
$p=8$, convergence with a factor of at least $0.1$ in each iteration.
Comparing the results for {\bf 2d-var} shown in Table~\ref{tab:2d-fan}
with the results for {\bf 2d-var$'$} in Table~\ref{tab:2d-fan2} shows
that the convergence does not degrade much for the
coefficient with 5-times smaller wavelength.

Next, we turn to the results for {\bf 3d-const} and {\bf 3d-var},
which we report in Tables~\ref{tab:3d-box} and \ref{tab:3d-fan},
respectively. For {\bf 3d-const}, all variants of the solver
converge. For this three-dimensional problem, the benefit of using
multigrid as preconditioner rather than as solver is even more evident
than in two dimensions.

Our results for {\bf 3d-var} are
summarized in Table~\ref{tab:3d-fan}. As for {\bf 2d-var}, the
performance of multigrid when used as a solver degrades
for orders $p>4$. We can also observe that the low-order matrix based on
the high-order node points represents a good preconditioner for the
high-order system. 

\begin{table}
  \caption{\label{tab:2d-fan} Iteration counts for two-dimensional
    warped-geometry, varying coefficient problem {\bf 2d-var} defined
    in \S\ref{subsec:tests}.  The finest
    mesh has $32\times 32$ elements and the multigrid hierarchy
    consists of three meshes.  For $p$-multigrid, the
    polynomial order is first reduced to $p=1$, followed by two
    geometric coarsenings of the mesh.
    For a
    detailed description of the different experiments reported in this
    table we refer to \S\ref{subsec:measures}.}  \centering
  \begin{tabular}{|r|c c|c c|c c||c c|c c|c c||c|} 
    \hline
    & \multicolumn{6}{c||}{MG as solver} & \multicolumn{6}{c||}{MG
      with pCG} & \!\!low-order MG\!\! \\
    \cline{2-13}
    \!\!\! order \!\!\!\! &  \multicolumn{2}{c|}{\!\!\scriptsize  Jacobi(3,3)\!\!} &  \multicolumn{2}{c|}{\!\!\scriptsize Cheb(3,3)\!\!} & \multicolumn{2}{c||}{\!\!\scriptsize  SSOR(2,1)\!\!} & \multicolumn{2}{c|}{\!\!\scriptsize Jacobi(3,3)\!\!} &  \multicolumn{2}{c|}{\!\!\scriptsize Cheb(3,3)\!\!} & \multicolumn{2}{c||}{\!\!\scriptsize SSOR(2,1)\!\!} & pCG\\
\hline
 & $h$ & $p$ & $h$ & $p$& $h$ & $p$& $h$ & $p$& $h$ & $p$& $h$ & $p$& ~ \\
 \cline{2-13}
1 & 14 & & 11 & & 6 & & 8 & & 7 & & 5 & & -   \\
2 & 20 & 19 & 15 & 15 & 7 & 8 & 10 & 10 & 8 & 8 & 5 & 6 & 16   \\
3 & 20 & & 16 & & 8 & & 10 & & 9 & & 6 & & 18  \\
4 & 22 & 21 & 21 & 19 & 10 & 9 & 11 & 10 & 10 & 10 & 7 & 6 & 19 \\
5 & -  & & 28 & & 12 & & 14 & & 12 & & 7 & & 21   \\
6 & -  & & 35 & & 13 & & 15 & & 13 & & 8 & & 23  \\
7 & -  & & 45 & & 16 & & 18 & & 15 & & 9 & & 24  \\
8 & -  & - & 52 & 46 & 17 & 15 & 20 & 20 & 16 & 15 & 9 & 8 & 25  \\
16 & - & - & 169 & 148 & 37 & 33 & 51 & 45 & 30 & 27 & 13 & 12 & 31 \\
\hline
  \end{tabular}
\end{table}

\begin{table}
  \caption{\label{tab:2d-fan2} Iteration counts for two-dimensional
    warped-geometry, varying coefficient problem {\bf 2d-var$'$}
    defined in \S\ref{subsec:tests}. This problem is identical to {\bf
      2d-var} (see Table~\ref{tab:2d-fan}),
    but the variations in the coefficient $\mu$ have a smaller wave length.
    The finest
    mesh has $32\times 32$ elements and the multigrid hierarchy
    consists of three meshes.  For $p$-multigrid, the
    polynomial order is first reduced to $p=1$, followed by two
    geometric coarsenings of the mesh.
    For a detailed description of the different experiments
    reported in this table we refer to \S\ref{subsec:measures}.} 
  \centering
  \begin{tabular}{|r|c c|c c|c c||c c|c c|c c||c|} 
    \hline
    & \multicolumn{6}{c||}{MG as solver} & \multicolumn{6}{c||}{MG
      with pCG} & \!\!low-order MG\!\! \\
    \cline{2-13}
    \!\!\! order \!\!\!\! &  \multicolumn{2}{c|}{\!\!\scriptsize  Jacobi(3,3)\!\!} &  \multicolumn{2}{c|}{\!\!\scriptsize Cheb(3,3)\!\!} & \multicolumn{2}{c||}{\!\!\scriptsize  SSOR(2,1)\!\!} & \multicolumn{2}{c|}{\!\!\scriptsize Jacobi(3,3)\!\!} &  \multicolumn{2}{c|}{\!\!\scriptsize Cheb(3,3)\!\!} & \multicolumn{2}{c||}{\!\!\scriptsize SSOR(2,1)\!\!} & pCG \\
\hline
 & $h$ & $p$ & $h$ & $p$& $h$ & $p$& $h$ & $p$& $h$ & $p$& $h$ & $p$& ~ \\
 \cline{2-13}
					1 & 14 & & 12 & & 8 & & 8 & & 8 & & 6 & & -  \\
					2 & 19 & 19 & 15 & 14 & 7 & 8 & 10 & 10 & 8 & 8 & 6 & 6 & 19 \\
					3 & 20 & & 17 & & 8 & & 10 & & 9 & & 6 & & 22  \\
					4 & 261 & 333 & 21 & 20 & 10 & 9 & 15 & 15 & 11 & 10 & 7 & 6 & 26  \\
					5 & - & & 30 & & 12 & & 19 & & 13 & & 8 & &  29  \\
					6 & - & & 39 & & 13 & & 37 & & 15 & & 8 & &  35  \\
					7 & - & & 52 & & 16 & & 78 & & 18 & & 9 & &  36  \\
					8 & - & - & 63 & 55   & 17 & 16 & 137 & 109 & 19 & 18 & 10 & 9  & 38 \\
				 16 & - & - & 232 & 201 & 67 & 76 &  -  &  -  & 44 & 37 & 19 & 18 & 56 \\
\hline
  \end{tabular}
\end{table}

\begin{table}
  \caption{\label{tab:3d-box} Iteration counts for three-dimensional unit cube
    problem  {\bf 3d-const} defined in \S\ref{subsec:tests}.
    The finest
    mesh has $8\times 8\times 8$ elements and the multigrid hierarchy
    consists of three meshes.  For $p$-multigrid, the
    polynomial order is first reduced to $p=1$, followed by two
    geometric coarsenings of the mesh.
    For a
    detailed description of the different experiments reported in this
    table we refer to \S\ref{subsec:measures}.}  \centering
	  \begin{tabular}{|r|c c|c c|c c||c c|c c|c c||c|} 
	    \hline
	    & \multicolumn{6}{c||}{MG as solver} &
            \multicolumn{6}{c||}{MG with pCG} &
            \!\!low-order MG\!\! \\
	    \cline{2-13}
	    \!\!\! order \!\!\!\! &  \multicolumn{2}{c|}{\!\!\scriptsize  Jacobi(3,3)\!\!} &  \multicolumn{2}{c|}{\!\!\scriptsize Cheb(3,3)\!\!} & \multicolumn{2}{c||}{\!\!\scriptsize  SSOR(2,1)\!\!} & \multicolumn{2}{c|}{\!\!\scriptsize Jacobi(3,3)\!\!} &  \multicolumn{2,1}{c|}{\!\!\scriptsize Cheb(3)\!\!} & \multicolumn{2}{c||}{\!\!\scriptsize SSOR(2,1)\!\!} & pCG \\
	\hline
	 & $h$ & $p$ & $h$ & $p$& $h$ & $p$& $h$ & $p$& $h$ & $p$& $h$ & $p$& ~\\
	 \cline{2-13}
1 & 6 & & 4 & & 4 & & 5 & & 4 & & 3 & & -   \\
2 & 8 & 8 & 4 & 5 & 4 & 5 & 6 & 6 & 4 & 4 & 4 & 4 &  25  \\
3 & 10 & & 7 & & 5 & & 6 & & 5 & & 5 & & 27  \\
4 & 11 & 10 & 8 & 7 & 6 & 5 & 7 & 7 & 6 & 5 & 5 & 4 & 28 \\
5 & 14 & & 10 & & 7 & & 8 & & 7 & & 5 & & 29  \\
6 & 16 & & 11 & & 7 & & 9 & & 7 & & 6 & & 32  \\
7 & 20 & & 15 & & 9 & & 10 & & 9 & & 6 & & 34 \\
8 & 22 & 19 & 17 & 15 & 9 & 8 & 10 & 10 & 9 & 8 & 6 & 6 & 35 \\
16 & 47 & 42 & 38 & 34 & 17 & 15 & 16 & 14 & 14 & 13 & 9 & 9 & 39 \\
\hline 
 \end{tabular}
\end{table}

\begin{table}
  \caption{\label{tab:3d-fan} Iteration counts for three-dimensional,
    warped-geometry, varying coefficient problem {\bf 3d-var} defined
    in \S\ref{subsec:tests}.     The finest
    mesh has $8\times 8\times 8$ elements and the multigrid hierarchy
    consists of three meshes.  For $p$-multigrid, the
    polynomial order is first reduced to $p=1$, followed by two
    geometric coarsenings of the mesh.
    For a detailed description of the
    different experiments reported in this table we refer to
    \S\ref{subsec:measures}.}  \centering
	  \begin{tabular}{|r|c c|c c|c c||c c|c c|c c||c|} 
	    \hline
	    & \multicolumn{6}{c||}{MG as solver} &
            \multicolumn{6}{c||}{MG with pCG} &
            \!\!low-order MG\!\! \\
	    \cline{2-13}
	    \!\!\! order \!\!\!\! &  \multicolumn{2}{c|}{\!\!\scriptsize  Jacobi(3,3)\!\!} &  \multicolumn{2}{c|}{\!\!\scriptsize Cheb(3,3)\!\!} & \multicolumn{2}{c||}{\!\!\scriptsize  SSOR(2,1)\!\!} & \multicolumn{2}{c|}{\!\!\scriptsize Jacobi(3,3)\!\!} &  \multicolumn{2}{c|}{\!\!\scriptsize Cheb(3,3)\!\!} & \multicolumn{2}{c||}{\!\!\scriptsize SSOR(2,1)\!\!} & pCG \\
	\hline
	 & $h$ & $p$ & $h$ & $p$& $h$ & $p$& $h$ & $p$& $h$ & $p$& $h$ & $p$& ~ \\
	 \cline{2-13}
1 & 13 & & 7 & & 5 & & 7 & & 5 & & 4 & & -  \\
2 & 17 & 18 & 13 & 13 & 7 & 7 & 9 & 9 & 8 & 8 & 5 & 5 & 26 \\
3 & 20 & & 16 & & 8 & & 10 & & 9 & & 6 & & 29  \\
4 & 23 & 22 & 18 & 18 & 9 & 9 & 11 & 11 & 9 & 9 & 7 & 6 & 31 \\
5 & 26 & & 21 & & 10 & & 12 & & 10 & & 7 & & 34  \\
6 & 30 & & 27 & & 12 & & 13 & & 12 & & 8 & & 37 \\
7 & 35 & & 34 & & 14 & & 14 & & 14 & & 8 & & 37  \\
8 & - & - & 40 & 38 & 16 & 15 & 18 & 17 & 15 & 14 & 9 & 9 & 38 \\
16 & - & - & 117 & 110 & 32 & 29 & 67 & 60 & 27 & 26 & 13 & 13 & 47\\
\hline
  \end{tabular}
\end{table}

\subsubsection{Mesh independence of iterations}\label{subsec:num_mesh}
To illustrate the mesh-independence of our multigrid-based solvers, we
compare the number of v-cycles required for the solution of the
two-dimensional problems {\bf 2d-const} and {\bf 2d-var} when
discretized on different meshes. In this comparison, the coarsest mesh
in the multigrid hierarchy is the same; thus, the number of levels in the
hierarchy increases as the problem is discretized on finer meshes. As
can be seen in Table~\ref{tab:meshInd}, once the mesh is sufficiently
fine, the number of iterations remains the same  for all polynomial orders.
\begin{table}[h]\centering
\caption{\label{tab:meshInd} Number of v-cycles required for the
  solution of the two-dimensional problems {\bf 2d-const} and {\bf
    2d-var} defined in \S\ref{subsec:tests} for different fine
  meshes and different polynomial orders. The coarsest grid for all
  cases has $2\times 2$ elements. In this comparison, multigrid with
  SSOR(2,1) smoothing is used as preconditioner in the conjugate
  gradient method. A star indicated that the
  corresponding test was not performed due to the large problem size.}
\begin{tabular}{|c|c|c|c|c|c|c|c||c|c|c|c|c|c|c|}
	\hline
	 & \multicolumn{7}{c||}{\bf 2d-const} &
        \multicolumn{7}{c|}{\bf 2d-var} \\  
  \hline
	order & 4  &  8 & 16 & 32 & 64 & 128 & 256 & 4  &  8 & 16 & 32 & 64 & 128 & 256 \\
	\hline 
	 1    &  3 &  4 &  4 &  4 &  4 &  4  &  4  & 3  &  4 & 5  &  5 &  5 &  5  &  5  \\
	 2    &  4 &  4 &  4 &  4 &  4 &  4  &  4  & 5  &  5 & 5  &  5 &  5 &  5  &  5  \\
	 4    &  5 &  5 &  4 &  4 &  4 &  4  &  4  & 6  &  6 & 7  &  7 &  7 &  7  &  7  \\
	 8    &  6 &  6 &  6 &  6 &  6 &  6  &  6  & 9  &  9 & 9  &  9 &  9 &  9  &  9  \\
	 16   &  9 &  9 &  9 &  9 &  9 &  *  &  *  & 13 & 13 & 13 & 13 & 13 &  *  &  *  \\
	 \hline
\end{tabular}

\end{table}

\subsubsection{Performance of block and $\ell_1$-Jacobi smoothers}\label{subsec:num_block}
For completeness, we also include a comparison with two common 
variants of the Jacobi smoother---the block-Jacobi and the
$\ell_1$-Jacobi point smoother.
We limit these comparisons to $8$ and $16$ order, and to the {\bf 2d-const}, {\bf
2d-var} and the {\bf 3d-var} problems. These results are summarized in Table~\ref{tab:block-jac}.

\paragraph{$\ell_1$-Jacobi smoother} These smoothers work by adding an appropriate diagonal matrix
to guarantee convergence \cite{BakerFalgoutKolevEtAl11}. They have the additional benefit of not requiring eigenvalue estimates 
compared with Chebyshev smoothers. In practice, while guaranteed convergence is desirable, the overall
work (i.e., number of iterations) increases.
In particular, point-Jacobi outperforms $\ell_1$-Jacobi as a smoother
for multigrid used as a solver as well as
a preconditioner for CG.

\paragraph{Block Jacobi smoother}
Schwarz-type domain decomposition smoothers are particularly promising
for high polynomial orders, such as order 8 or higher. Results
obtained with an elementwise block Jacobi preconditioner for orders 8
and 16 are summarized in Table~\ref{tab:block-jac}. For this
comparison, we invert the element matrices exactly, which can be
problematic with respect to computational time as well as storage for
realistic problems, in particular for warped meshes and high
polynomial orders. One remedy is to use approximate inverse element
matrices \cite{LottesFischer05}. As can be seen in
Table~\ref{tab:block-jac}, the
number of iterations is reduced compared to pointwise Jacobi
smoothing; however, this does not imply a faster method since
block-Jacobi smoothing is, in general, more expensive. Again, a
high-performance implementation is required to assess the
effectiveness of the different methods.
\begin{table}\centering
	\caption{\label{tab:block-jac} Comparison between different Jacobi
          smoothers---point, block and $\ell_1$. Shown is the number of
          iterations, obtained with 3
          pre- and 3 post-smoothing steps. All experiments used a damping factor of $\omega=2/3$.}
\begin{tabular}{|c|c|c|c|c|c|c|c|c|c|c|c|c|c|c|c|c|c|c|}
	\hline
	\!\!order\!\! & \multicolumn{6}{c|}{\bf 2d-const} &
        \multicolumn{6}{c|}{\bf 2d-var} & \multicolumn{6}{c|}{\bf 3d-var} \\
	\cline{2-19}
  & \multicolumn{3}{c|}{MG} & \multicolumn{3}{c|}{pCG} &\multicolumn{3}{c|}{MG} & \multicolumn{3}{c|}{pCG}  & \multicolumn{3}{c|}{MG} & \multicolumn{3}{c|}{pCG}  \\
	\hline
	 & {\!pt\!} & {\!blk\!} & {\!$\ell_1$\!} & {\!pt\!} & {\!blk\!} & {\!$\ell_1$\!} & {\!pt\!} & {\!blk\!} & {\!$\ell_1$\!} & {\!pt\!} & {\!blk\!} & {\!$\ell_1$\!}& {\!pt\!} & {\!blk\!} & {\!$\ell_1$\!}& {\!pt\!} & {\!blk\!} & {\!$\ell_1$\!}\\
  \hline \small
  8  & 17 & 16 & 51 & 9  &  8 & 16 & - & 31 & \!111\! & - & 12 & 57 & - & 30 & \!176\! & 18 & 13 & 37 \\    
  16 & 40 & 31 & \!133\! & 14 & 12 & 27 & - & 61 & - & 51 & 17 & \!186\! & - & 52 & 48 & 67 & 17 & 68  \\   
  \hline
\end{tabular}
\end{table}
In the next section, we summarize our findings and draw conclusions.

\section{Discussion and conclusions}
\label{sec:discuss}

Using multigrid as preconditioner in the conjugate gradient (CG)
method rather than directly as solver results in significantly faster
convergence, which more than compensates for the additional work
required by the Krylov method.  This is particularly true for
high-order methods, where the residual computation is more expensive
than for low-order methods, thus making the additional vector
additions and inner products in CG negligible.  For problems with
varying coefficients, we find that the number of v-cycles
decreases by up to a factor of three when multigrid is combined with
the conjugate gradient method.

None of the tested approaches yields a number of iterations that is
independent of the polynomial order; Nevertheless, point smoothers can
be efficient for finite element discretizations with polynomial orders
up to $p=16$.  For constant coefficient, all tested multigrid
hierarchy/smoother combinations (Jacobi, Chebyshev-accelerated Jacobi
and Gauss-Seidel SSOR smoothing) lead to converging multigrid
methods. In general, the difference in the number of iterations
between $h$- and $p$-multigrid is small.
Problems with strongly varying coefficients on deformed
geometries are much more challenging. Here, SSOR outperforms
Jacobi-based smoothers for orders $p>4$. However, in a distributed
environment, where Gauss-Seidel smoothing is usually more difficult to
implement and requires more parallel communication,
Chebyshev-accelerated Jacobi smoothing represents an interesting
alternative to SSOR. It is as simple to implement as Jacobi smoothing
but requires significantly less iterations to converge; compared to
point Jacobi smoothing, it additionally only requires an estimate of
the largest eigenvalue of the diagonally preconditioned system matrix.

We find that a low-order operator based on the high-order node points
is a good preconditioner, and it is particularly attractive for
high-order discretizations on unstructured meshes, as also observed in
\cite{Brown10, DevilleMund90, HeysManteuffelMcCormickEtAl05}.  When
combined with algebraic multigrid for the low-order operator, the
smoother on the finest mesh can either use the low-order or the
high-order residual. Initial numerical tests indicate that the latter
choice is advantageous, but this should be studied more systematically.

\section*{Acknowledgments}
We would like to thank Tobin Isaac for helpful discussions on the
low-order preconditioner. Support for this work was
  provided through the U.S.~National Science Foundation (NSF) grants
  CMMI-1028889 and   
  ARC-0941678,       
 and through the Scientific Discovery through Advanced
  Computing (SciDAC) projects 
  DE-SC0009286,   
  and DE-SC0002710 
  funded by the U.S.~Department of Energy
  Office of Science, Advanced Scientific Computing Research and
  Biological and Environmental Research.

\bibliographystyle{siam}
\bibliography{ccgo,mg}

\end{document}